\RequirePackage{rotating}

\documentclass{amsproc}

\usepackage[english]{babel}

\PassOptionsToPackage{pagebackref}{hyperref}

\usepackage{def-packages}
\usepackage{def-makros}
\usepackage{def-thm}

\begin{document}


\title{Culler-Shalen seminorms of fillings of the Whitehead link exterior}
\author{Gabriel Indurskis}
\address{D\'epartement de math\'ematiques\\
201, Avenue du Pr\'esident-Kennedy\\
Montr\'eal (Qu\'ebec)\\ Canada H2X 3Y7}
\curraddr{Department of Mathematics\\
 Champlain College St-Lambert\\
 900 Riverside Drive\\
 St-Lambert (Qu\'ebec)\\ Canada J4P 3P2}
\email{gindurskis@crcmail.net} 
\thanks{This work was supported in part by a
  scholarship of the \emph{Institut des Sciences Math\'ematiques},
  Montr\'eal.}

\begin{abstract}
  We determine the total Culler-Shalen seminorms for the $3$-manifolds
  $\Wpq:=W(p/q,\ \cdot\ )$ obtained by Dehn filling with slope $p/q$ on one
  boundary component of the Whitehead link exterior $W$ when $p$ is odd.  As
  part of the proof, we use an explicit parametrization of the eigenvalue
  variety of $W$ to find a one-variable polynomial whose roots characterize
  characters of p-reps of $\pi_{1}(\Wpq)$, i.e.\ representations with values
  in $\SL_{2}(\C)$ which are parabolic on the peripheral subgroup.
\end{abstract}
\subjclass{Primary 57M27; Secondary 20C15, 57N10, 57R65}
\keywords{3-manifold, Dehn surgery, character variety, Culler-Shalen seminorm}
\maketitle




\section{Introduction}
\label{sec:introduction}

\begin{figure}[htbp]
  \centering
  \input{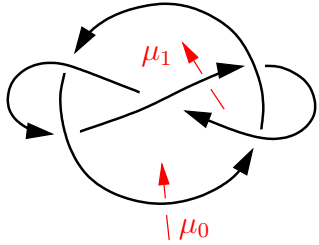}
  \\[-3ex]\caption{The right-handed Whitehead link}
  \label{fig:whitehead-link-right}
\end{figure}

Let $W$ be the exterior of the right-handed Whitehead link in $S^{3}$ with the
canonical framing $(\mu_{0}, \lambda_{0})$ and $(\mu_{1}, \lambda_{1})$ as
indicated in figure~\ref{fig:whitehead-link-right}.  Let
$\Wpq:=W(\pqslope,\ \cdot\ )$ be the manifold obtained from $W$ by Dehn
filling the first boundary torus along the slope $\pqslope$, i.e.\ by
attaching a solid torus by identifying its meridian with the curve
$\mu_{0}^{p}\lambda_{0}^{q}$.  Here $p$ and $q$ are coprime integers, and
unless otherwise noted we will fix the convention $q>0$. (Note that there is
an ambient isotopy taking one component of the link to the other, so in fact
$W(\pqslope,\ \cdot\ )\isom W(\ \cdot\ ,\pqslope)$. To avoid confusion, we
will nevertheless always consider the boundary with curves $\mu_{0}$,
$\lambda_{0}$ to be filled, and the curves $\mu_{1}$ and $\lambda_{1}$ to lie
on the remaining boundary torus of $\Wpq$ after filling.)

The manifolds $\Wpq$ constitute a rich class of $3$-manifolds whose boundary
consists of one torus, and contain important infinite subfamilies of
well-known manifolds such as the twist knot exteriors (for
$p/q=1/n$, $n\in\Z$) and an important class of once-punctured torus
bundles (for $p/q=m\in\Z$, cf.\
\cite{MR93i:57019}). Note that for example the figure-$8$
knot complement ($p/q=-1$) and its so-called
sister manifold ($p/q=5$) are part of this family. Since the Whitehead link
exterior admits a complete hyperbolic structure of finite volume, the same is
true for all but a finite number of the manifolds $\Wpq$ (cf.\
\cite{ThurstonNotes}). Furthermore, unless $p/q\in\{0,4\}$, all manifolds
$\Wpq$ are \emph{small} manifolds,
i.e.\ they do not contain any closed essential surfaces (cf.\
\cite[Appendix]{MR1913944}).

The goal of this paper is to determine the total Culler-Shalen seminorm of the
filled manifold $\Wpq$ when $p$ is odd, in terms of the filling slope $p/q$.
This work can be seen as a generalization of similar work previously done for
twist knots by Boyer, Mattman and Zhang (\cite{MR1664959}).  The results given
in this paper completely include these previous results, and extend them to a
much larger class of fillings.

\begin{Theorem}\label{thm:total-cs-norm-Wpq}
  Let $p$ and $q$ be coprime integers with $q>0$ and $p$ odd,\
$p/q\neq 3$.
 The total Culler-Shalen seminorm of the manifold
  $\Wpq=W(\pqslope,\ \cdot\ )$ is given by%
  \begin{align*}
    \CSnorm{\gamma} &=  \sum_{j=1}^{3} a_{j}\dist(\gamma,\beta_{j}), 
  \end{align*}
  where $\dist(\cdot,\cdot)$ denotes the geometric intersection number of two
  slopes, the slopes $\beta_{j}$ ($j=1,2,3$) are the possible boundary slopes
  for $\Wpq$ given in table~\ref{tab:boundary-slopes-Wrh-summary} below,
  and the coefficients
  $a_{1}$, $a_{2}$, and $a_{3}$ are given in
  table~\ref{tab:summary-cs-norm-p-odd}. Furthermore, the minimal non-zero
  total norm $s$ is as given in the last column of
  table~\ref{tab:summary-cs-norm-p-odd} on page~\pageref{tab:summary-cs-norm-p-odd}.
\end{Theorem}
We remark that the condition that $p$ be odd (which is required in the proofs
of section~\ref{sec:whitehead-seifert-fillings} to avoid a possible
obstruction to lifting $\PSL_{2}(\C)$ representations) automatically excludes
the special cases $W_{0}$ and $W_{4}$ (both \emph{large}, toroidal manifolds
which do not admit a complete hyperbolic structure), as well as $W_{6}$:
Filling this manifold a second time with slope $1$ produces the closed
manifold $L_{2}\connsum L_{3}$, a connected sum of lens spaces, causing the
arguments in section~\ref{sec:whitehead-seifert-fillings} to fail. The same
argument applies to $W_{3}$, which therefore has to be excluded separately.

\begin{table}[htb]
  \centering
  \begin{tabular}{C||C|CCC|C}
   r=\pqslope\in  & \beta_{1} & \beta_{2} &&&  \beta_{3} \\
\hline
 & & &&& \\[-2ex]
{[-\infty,0]}& 4 & \ds 4r^{-1}=\frac{4q}{p}       &\in&{[-\infty,0]} & 0\\[3ex]
{[0,2]}      & 4 & \ds 4r^{-1}+2=\frac{2p+4q}{p}&\in&{[4,\infty]} & 0\\[3ex]
{[2,4]}      & 4 & \ds 6-r=\frac{-p+6q}{q}&\in&{[2,4]} & 0\\[3ex]
{[4,\infty]} & 4 & \ds 4(r-2)^{-1}=\frac{4q}{p-2q}&\in&{[0,2]} & 0\\
  \end{tabular}
  \\[2ex]\caption{Possible boundary slopes of $\Wpq$}
  \label{tab:boundary-slopes-Wrh-summary}
\end{table}

\begin{table}[htb]
  \centering
  \begin{tabular}{C||C|C|C|C}
  \pqslope\in & a_{1} & a_{2} & a_{3} & s\\
  \hline
  (-\infty,0)               &  -p+2q-1  & 2  &  2q-2  & -3p+4q-3 \\
  (0,2)                     &  -p+2q-1  & 2  &  2q-2  &   p+4q-3 \\
  (2,3)\cup (3,4)           &   p-2q-1  & 4  &  2q-2  &   p+4q-3 \\
  (4,6)\cup (6,\infty)      &   p-2q-1  & 2  &  2q-2  &  3p-4q-3 \\
  \end{tabular}
  \\[2ex]\caption{The parameters for the total Culler-Shalen seminorm of $\Wpq$ for $p$ odd, $q>0$}
  \label{tab:summary-cs-norm-p-odd}\vspace*{-.6cm}
\end{table}

The knowledge of the coefficients $a_{j}$ given in
Theorem~\ref{thm:total-cs-norm-Wpq} has immediate consequences: If $a_{j}>0$,
the fundamental polygon (the boundary of the norm-ball of radius $s$) of the
Culler-Shalen seminorm has a vertex corresponding to $\beta_{j}$, and
Lemma~6.1 of \cite{BZfinite} then shows that in this case $\beta_{j}$ is a
\emph{strongly detected}, strict boundary slope, i.e.\ the slope of an
essential surface associated to some ideal point of the character variety of
$\Wpq$, which is neither a fibre nor semi-fibre. Inspection of
table~\ref{tab:summary-cs-norm-p-odd} therefore implies:%
\begin{Corollary}\label{cor:Wpq-boundary-slopes-are-detected}%
  Let $p$ and $q$ be coprime integers with $q>0$, $p$ odd, and $p/q\neq 3$.
\begin{enumerate}
\item $\beta_{1}=4$ is a strongly detected slope of $\Wpq$ \emph{unless} $p=2q\pm 1$.
\item $\beta_{2}$ is \emph{always} a strongly detected slope of $\Wpq$.
\item $\beta_{3}=0$ is a strongly detected slope of $\Wpq$ iff $q>1$, otherwise it
  is the slope of a fibre in the once-punctured torus bundle $W(p)$.
\end{enumerate}
\end{Corollary}

An interesting by-product of the proof of Theorem~\ref{thm:total-cs-norm-Wpq}
is the complete determination of the so-called p-reps of the fundamental group
of $\Wpq$, which are of interest in their own right. These are representations
into $\SL_{2}(\C)$ which are parabolic on the peripheral subgroup, and were
first studied for knot groups by Riley (cf.\ \cite{MR0300267}). (To be
precise, the generator corresponding to the meridian of the boundary torus is
sent to a non-trivial element of trace $\pm 2$ in $\SL_{2}(\C)$. As a
consequence, the generator corresponding to the longitude will be either
parabolic as well, or $\pm I$.) In his classical work, he showed how under
certain conditions on the group presentation (e.g.\ for the fundamental group
of a $2$-bridge knot) all such p-reps can be characterized by the complex
roots of a single one-variable polynomial. This is in contrast to the general
situation, where only a characterization by a system of \emph{several}
polynomial equations can be expected.

Passing through the so-called \emph{eigenvalue variety}
 (as defined by
Tillmann in \cite{MR2104008}) of the unfilled Whitehead link exterior $W$\!, we are able
to find a similar approach which seems naturally adapted to the situation of
Dehn filling to determine a \emph{single} one-variable polynomial whose roots are in
two-to-one correspondence with the set of eigenvalues of p-reps. This in
turn is shown to be in bijective correspondence with the set of p-reps
itself, and we obtain a complete characterization of all p-reps of $\Wpq$:
\begin{Theorem}\label{thm:conjugacy-classes-preps-charact-by-roots-of-respq}
  The conjugacy classes of irreducible p-reps of the manifold $\Wpq$ ($p$ can
  be odd or even) are in two-to-one correspondence with the non-trivial roots
  $s\notin \{0, \pm 1\}$ of the Laurent polynomial
\[ \res_{p,q}(s)= s^{p-2q} + (-1)^{q+1}2T_{q}(y) + s^{-p+2q}, \] where
$y=\tfrac{1}{2}(-s^{2}+4-s^{-2})$, and $T_{q}(y)$ is the $q$-th Chebyshev polynomial of first
type. The conjugacy classes of non-abelian reducible p-reps are in two-to-one
correspondence with the non-trivial roots of unity $s\neq \pm 1$, $s^{p}=1$.
\end{Theorem}

We remark that whereas it seems a hard analytic problem to determine if the
non-trivial roots of these polynomials are always simple, this is in fact a
consequence of Theorem~\ref{thm:total-cs-norm-Wpq} when $p$ is odd:
\begin{Corollary}\label{cor:roots-of-res-are-simple}
  When $p$ is odd (including $p/q=3$), all non-trivial roots $s\notin \{0,\pm 1\}$ of $\res_{p,q}$ are simple.
\end{Corollary}
\begin{proof}
Note that even though Theorem~\ref{thm:total-cs-norm-Wpq} excludes $p/q=3$, we
can include it here simply by manually inspecting $\res_{3,1}(s)=-s^{3}(s+1)^{2}(s^{2}-3s+1)$.

For all other cases with $p$ odd, Theorem~\ref{thm:total-cs-norm-Wpq} shows
that the upper bounds for the count of distinct non-trivial roots of
$\res_{p,q}$ obtained in section~\ref{sec:roots-resultant} are actually
attained, which can only happen if these roots are all simple.
\end{proof}

We note that experimental evidence seems to \emph{suggest} that this also
holds when $p$ is even. Furthermore, we would like to point out that all our
results would immediately extend to the case when $p$ is even if a work-around
for the lifting obstruction in section~\ref{sec:whitehead-seifert-fillings}
could be found.

\bigskip
The organization of this paper is as follows:
\begin{description}
\item[Section~\ref{sec:outline}] After discussing some preliminaries about
  character varieties and the Culler-Shalen seminorm, we give the
  general outline for the proof of Theorem~\ref{thm:total-cs-norm-Wpq}. We
  also determine the set of possible boundary slopes of $\Wpq$.

\item[Section~\ref{sec:finding-preps}] In this section, we identify the
  p-reps of the manifold $\Wpq$ with points in the eigenvalue variety of the
  unfilled manifold $W$ characterized by the roots of a single one-variable
  Laurent polynomial $\res_{p,q}$, proving
  Theorem~\ref{thm:conjugacy-classes-preps-charact-by-roots-of-respq}.  A
  detailed analysis of the properties and symmetries of these recursively
  defined polynomials is given, including an upper bound for the number of
  their roots in terms of $p$ and $q$, which in turn gives an upper bound for
  the number of conjugacy classes of p-reps.

\item[Section~\ref{sec:minimal-CS-norm}] We show that the \emph{minimal}
  nonzero value of the Culler-Shalen seminorm of $\Wpq$ is precisely equal to
  the number of conjugacy classes of p-reps, and therefore satisfies the same
  bounds as shown in the previous section: Using group cohomology calculations
  and the characterization of p-reps from section~\ref{sec:finding-preps}, it
  is shown that each character of a p-rep is a smooth point on a non-trivial
  curve in the character variety of $\Wpq$, and that it furthermore is a
  simple root of the trace function $f_{\mu_{1}}$.

\item[Section~\ref{sec:whitehead-seifert-fillings}] We explicitly determine
  the Culler-Shalen seminorm of slopes $\sigma$ resulting in Seifert-fibered
  manifolds $\Wpq(\sigma)$ when $p$ is odd. This is done by a detailed study of
  $\PSL_{2}(\C)$ representations of the orbifold fundamental groups of their base
  orbifolds, which can always be lifted to $\SL_{2}(\C)$ representations when
  $p$ is odd, but not necessarily if $p$ was even.

\item[Section~\ref{sec:total-norm-solving-systems}] The total Culler-Shalen
  seminorm is determined for $p$ odd: Combining the results from the previous
  sections, we solve the linear systems relating the seminorms obtained for
  the Seifert fillings with the coefficients determining the seminorm.

\end{description}

\subsection*{Acknowledgments}
\label{sec:acknowledgments}
\vspace*{-1ex}
{\small
  \emph{The work contained in this paper forms part of the author's dissertation
  and could not have been achieved without the never-waning support of his
  supervisor, Steven Boyer, for which the author would therefore like to express his sincerest
  thanks. Further thanks go to Stephan Tillmann for suggesting the
  ``detour'' through the eigenvalue variety and for many interesting
  discussions. Finally, the author would like to thank the eagle-eyed
  anonymous reviewer for catching a potentially show-stopping sign mistake,
  and whose comments and recommendations
  have greatly improved the readability of this paper.  
}


\section{Preliminaries and outline}
\label{sec:outline}

We briefly recall the construction of the Culler-Shalen
norm developed in \cite{CGLS},
 and its generalization to a seminorm by
 Boyer and Zhang in
\cite{BZseminorms}.

Let $M$ be a connected, orientable, compact, irreducible and boundary-irredu\-cible
$3$-manifold with boundary $\boundary M$ consisting of one
torus. Unless noted otherwise, we will generally suppress the base point of
fundamental groups.

Denote by $\RV(M):=\Hom(\pi_{1}(M); \SL_{2}(\C))$ the
\emph{representation variety} of $M$\!, and by $\CV(M)$ its \emph{character variety},
i.e.\ the set of characters of representations in $\RV(M)$ (cf.\
\cite{MR683804}). These two sets are in fact affine algebraic varieties and as
such can be identified with subsets of $\C^{N}$ (for appropriate $N\in \N$),
defined by a set of polynomial equations. We will occasionally need to
consider representations into $\PSL_{2}(\C)$, which we will denote by 
$\prho\in \PRV(M):=\Hom(\pi_{1}(M); \PSL_{2}(\C))$.

Since $M$ is boundary-irreducible, we can identify $\pi_{1}(\boundary M)$ with its image
under the natural inclusion $\pi_{1}(\boundary M)\ra\pi_{1}(M)$ (defined up to
conjugation). Denote the first homology of the boundary with integer coefficients
$H_{1}(\boundary M; \Z)$  by $L$ and regard it as
a lattice in the $2$-dimensional real vector space $V=H_{1}(\boundary M;
\R)$. Using the Hurewicz homomorphism, we can identify $L$ with
$\pi_{1}(\boundary M)$.

A \emph{slope} $r$ is the isotopy class of an essential, unoriented, simple
closed curve on $\boundary M$. It can be identified with a pair of elements
$\pm \delta=\pm(p,q)\in L$, or with the rational number $p/q\in
\Q\union\{\infty\}$, where $\infty=1/0$ and $-\infty=-1/0$ are identified.
Let $\alpha(r)$ be either one of $\pm \delta$. We will often slightly abuse
notation and interchange freely between these notations; for example writing
$M(\delta)$ or $M(p/q)$ instead of $M(r)$ to denote the manifold obtained by
Dehn filling $M$ along the slope $r$.  A slope is called a \emph{boundary
  slope} if it represents the isotopy class of the boundary components (on one
component) of $\boundary M$ of an essential surface $S\subset M$. A boundary
slope is called \emph{strict} if it is the
boundary slope of an essential surface $S$ which is neither a fibre or a
semi-fibre in a fibre or semi-fibre bundle structure on $M$. For two elements
$\delta_{1},\delta_{2}\in L$, we denote their geometric intersection number,
called their \emph{distance}, by $\dist(\delta_{1},\delta_{2})$. If
$\delta_{i}=(p_{i},q_{i})$, we have $\dist(\delta_{1},\delta_{2})=|p_{1}q_{2}-p_{2}q_{1}|$.

Let $\CV_{0}\subset\CV$ be a curve in the character variety $\CV=\CV(M)$,
i.e.\ an irreducible algebraic component of complex dimension $1$.  Call
$\CV_{0}$ \emph{non-trivial} if it contains the character of an irreducible
representation.

For $\gamma\in\pi_{1}(M)$, define the regular function $I_{\gamma}\map
\CV_{0}\ra\C$ by $I_{\gamma}(\chi_{\rho})=\chi_{\rho}(\gamma)=\trace
\rho(\gamma)$. Use this to define the map $f_{\gamma}:=I_{\gamma}^{2}-4$,
which is also regular and can hence be pulled back to the smooth projective
model (i.e.\ the unique non-singular projective variety birationally
  equivalent to the curve $\CV_{0}$, cf.\ \cite[Ch.~II, \S~4.5 and
  5.3]{MR1328833})  $\CVsp_{0}$ of $\CV_{0}$ (we will denote this map also by
$f_{\gamma}$).

For $\gamma\in L$, define 
$\CSnorm{\gamma}_{\CV_{0}}$  as the degree of
the map $f_{\gamma}\map \CVsp_{0}\ra \CP^{1}$. This can be extended
naturally to define a function on the real vector space $V:=H_{1}(\boundary M;\R)$:
\begin{equation}
  \CSnorm{\cdot}_{\CV_{0}}\map H_{1}(\boundary M;\R)\ra [0;\infty)
\end{equation}

This function is a seminorm on $H_{1}(\boundary M;\R)$, called the
\emph{Culler-Shalen seminorm} associated to the curve $\CV_{0}$ (cf.\
\cite[Proposition~1.1.2]{CGLS} for the case of a norm, and
\cite[Section~5]{BZseminorms} for the general case).  Seen as a seminorm on a
$2$-dimensional real vector space, $\CSnorm{\cdot}_{\CV_{0}}$ is either a
norm, a nonzero indefinite seminorm, or identically zero. When $M$ is small,
the last possibility cannot occur (cf.\ \cite[Proposition~5.5]{BZseminorms}).
Furthermore $\dim \CV(M)\leq 1$ in this case (cf.\ \cite{CCGLS}), and it is
therefore natural to consider the sum of the Culler-Shalen seminorms over all
\emph{non-trivial} curves $\CV_{i}\subset\CV(M)$, which we will call the
\emph{total} Culler-Shalen seminorm:
\begin{equation}
  \CSnorm{\gamma}:=\sum_{i} \CSnorm{\gamma}_{i}
\end{equation}
We denote the union of all non-trivial curves in $\CV(M)$ by $\CVi(M)$.

If $f_{\gamma}$ is non-constant on $\CV_{0}$ for some $\gamma\in L$, define
\begin{equation}
  s_{\CV_{0}}:=  \min \{ \CSnorm{\gamma}_{\CV_{0}} \where \gamma\in L,\
      \gamma\neq 0,\ 
      \text{$f_{\gamma}$ not constant on $\CV_{0}$} \},
\end{equation}
otherwise set $s_{\CV_{0}}:=0$. We call $s_{\CV_{0}}$ the \emph{minimal}
seminorm on $\CV_{0}$, the corresponding value for the total seminorm is the
\emph{total minimal seminorm}, denoted by $s$.

The following lemma enables us to
determine the Culler-Shalen seminorm from only a few known values:%
\begin{Lemma}[{\cite[Lemma~6.2]{BZfinite}}]\label{lem:CSnorm-as-lin-comb}
  If $\beta_{j}$ are the boundary slopes of $M$ and $\CV_{i}\subset\CV(M)$ is
  a non-trivial curve (i.e.\ a curve which contains the character of an
  irreducible representation), then \[\CSnorm{\gamma}_{\CV_{i}} = \sum_{j}
  a_{ij}\dist(\gamma,\beta_{j}),\] where the $a_{ij}$ are non-negative, even
  integers independent of $\gamma\in \pi_{1}(M)$.
\qed
\end{Lemma}
For the total seminorm, we therefore get
\begin{equation}\label{eq:total-CSnorm-as-linear-comb-of-dist}
  \CSnorm{\gamma} = \sum_{j} a_{j}\dist(\gamma,\beta_{j}),
\end{equation}
with non-negative even integers $a_{j}$.  

To determine the total seminorm for the manifolds $\Wpq$, it therefore
suffices to determine the boundary slopes $\beta_{j}$ and the corresponding
coefficients $a_{j}$. Following a technique due to Mattman (cf.\
\cite{MR1949779}), we will use explicit calculations of the seminorm of
certain slopes (cf.\ section~\ref{sec:whitehead-seifert-fillings}) together with
equation~\eqref{eq:total-CSnorm-as-linear-comb-of-dist} to set up a system of
linear equations with unknowns $a_{j}$ and $s$. In our case, this system turns
out to be underdetermined, requiring one more piece of information to
successfully solve it. We therefore first study the set of p-reps of $\Wpq$ (cf.\
section~\ref{sec:finding-preps}) and use
this knowledge to explicitly determine the minimal total seminorm $s$ (cf.\
section~\ref{sec:minimal-CS-norm}). With this information, we are then able to
solve the linear systems in section~\ref{sec:total-norm-solving-systems},
completing the proof of Theorem~\ref{thm:total-cs-norm-Wpq}.

We finish this section by determining the set of possible boundary slopes of $\Wpq$:

\subsection{Boundary slopes of the unfilled Whitehead link exterior}
\label{sec:bound-slop-unfilled}

In \cite{MR2299787}, Hoste and Shanahan give an explicit
algorithm, building on work of Floyd and Hatcher (cf.\
  \cite{MR924770}) and on the unpublished dissertation of Lash (cf.\
  \cite{L1993}), to determine the boundary slopes of $2$-bridge links. In
particular, they show that up to symmetry between the two boundary tori, the
\emph{left-handed} Whitehead link ($L_{3/8}$ in $2$-bridge notation) has
boundary slope pairs $(0,0)$, $(0,\emptyset)$, $(-4,-2)$, $(-4,\emptyset)$
(here $\emptyset$ signifies that the intersection with one of the boundary
components is empty) and the continuous spectrum of slope pairs
$(2t^{-1},2t)$, $(-2t^{-1},-2-2t)$, $(-3+s,-3-s)$ parametrized by $t\in
(\Q\union \{\infty\})\intersect [0,\infty]$ and $s\in \Q\intersect[-1,1]$
(cf.\ table~4 in \cite{MR2299787}, using $k=1$).

After applying an orientation-reversing
homeomorphism to switch to the \emph{right-handed} Whitehead link (which
simultaneously changes the signs of all slopes), we use the symmetry between
the two boundary components and reparametrize the slope pairs such that the
first is always equal to $r=p/q$. The results are shown in
table~\ref{tab:boundary-slopes-Wrh}.

\begin{table}[htbp]
  \centering
  \begin{tabular}{C|CCC}
 \text{$r=p/q$ on $T_{0}$}   & \text{$p'/q'$ on $T_{1}$} &&  \\\hline\hline
 0                      &  0\text{ or }\emptyset &&\\
 4                      &  2\text{ or }\emptyset &&\\\hline
 \emptyset  &  0 &&\\
 2\text{ or }\emptyset  &  4 &&\\\hline
&&&\\[-1ex]
\ds r  \in {[-\infty,0]}& \ds4r^{-1}=\frac{4q}{p}       &\in&{[-\infty,0]}\\[3ex]
\ds r  \in {[0,2]}      & \ds4r^{-1}+2=\frac{2p+4q}{p}&\in&{[4,\infty]}\\[3ex]
\ds r  \in {[2,4]}      & \ds6-r=\frac{-p+6q}{q}&\in&{[2,4]} \\[3ex]
\ds r  \in {[4,\infty]} & \ds4(r-2)^{-1}=\frac{4q}{p-2q}&\in&{[0,2]} \\
  \end{tabular}
  \\[2ex]\caption{Boundary slope pairs of the right-handed Whitehead link $W$}
  \label{tab:boundary-slopes-Wrh}
\end{table}

\subsection{Boundary slopes of the once-filled Whitehead link exterior}
\label{sec:boundary-slopes-Wpq}

Denote the two boundary components of $W$ by $T_{0}$ and $T_{1}$.
We now proceed to determine the boundary slopes of the once-filled Whitehead
link exterior $\Wpq=W(\pqslope,\ \cdot\ )=W\union_{\mu_{0}^{p}\lambda_{0}^{q}}
V$, where $V$ is the attached solid torus. Note that we can view $W$ as lying
inside $\Wpq$.
\begin{Lemma}
  Let $S$ be an essential surface in $\Wpq$. Then either $S$ is in fact an
  essential surface in $W\subset \Wpq$, or it is the cap-off of an
  essential surface $S'$ in $W\subset \Wpq$. In particular, $S$ has non-empty
  intersection with at least one of $T_{0}$ and $T_{1}$.
\end{Lemma}
\begin{proof}
  Let $S$ be an essential surface in $\Wpq$.  We can isotope $S$ so it
  transversely intersects $T_{0}\subset W\subset \Wpq$ and $T_{1}\subset \Wpq$,
  as well as the core of the filling torus $V$.  If $S\intersect
  T_{0}=S\intersect T_{1}=\emptyset$, $S$ can be viewed as a closed essential
  surface in $W$, which is a contradiction, since $W$ is small.  We therefore
  see that $S$ has to intersect at least one of $T_{0}$ or $T_{1}$.
\end{proof}

We now have two cases to consider:
\begin{description}
\item[$S\intersect \interior V=\emptyset$] In this case, $S$ is in fact an
essential surface in $W\subset\Wpq$ with empty intersection with $T_{0}$.
According to table~\ref{tab:boundary-slopes-Wrh}, it therefore must have
boundary slope on $T_{1}$ equal to either $0$ or $4$.

\item[$S\intersect\interior V\neq \emptyset$] In this case, the intersection
  $S\intersect T_{0}$ must consist of simple closed curves of slope
  $\pqslope$. This implies that $S'=S\intersect W$ is an essential surface in
  $W$ with boundary slope $\pqslope$ on $T_{0}$. Unless $\pqslope\in \{0,4\}$,
  this implies that $S'\intersect T_{1}\neq \emptyset$, and we get the
  corresponding boundary slope on $T_{1}$ from
  table~\ref{tab:boundary-slopes-Wrh}, dependent on $\pqslope$.  We remark
  that if $\pqslope\in \{0,4\}$, $S$ could in fact be a closed essential torus
  in $\Wpq$.
\end{description}

We therefore have at most three \emph{candidates} for boundary slopes in the
manifold $\Wpq$, as summarized in table~\ref{tab:boundary-slopes-Wrh-summary}
on page~\pageref{tab:boundary-slopes-Wrh-summary}: We have two fixed slopes
$\beta_{1}=4$ and $\beta_{3}=0$, and one varying slope $\beta_{2}$ which
depends on $p$ and $q$. The parametrization of the latter depends on the range
of the filling slope $\pqslope$. 

Note that we have at this point not determined if these slopes are in fact
boundary slopes of $\Wpq$. For the following, we will rather only assume 
that the set of boundary slopes of $\Wpq$ is \emph{contained} in
the set $\{\beta_{1},\ \beta_{2},\ \beta_{3}\}$.  The proof of
Theorem~\ref{thm:total-cs-norm-Wpq} will then in turn imply
Corollary~\ref{cor:Wpq-boundary-slopes-are-detected}, showing that unless
$p=2q\pm 1$ or $q=1$, all three slopes are indeed strongly detected, strict
boundary slopes.




\section{The set of p-reps of a once-filled Whitehead link exterior}
\label{sec:finding-preps}


The goal of this section is to prove
Theorem~\ref{thm:conjugacy-classes-preps-charact-by-roots-of-respq} (in
sections~\ref{sec:detour-through-eigenv}--\ref{sec:irred-solut}), as well as
to provide upper bounds for the number of non-trivial roots of the Laurent
polynomials $\res_{p,q}$ (in section~\ref{sec:roots-resultant}), which in turn
then gives upper bounds for the number of conjugacy classes of p-reps (in
section~\ref{sec:number-p-reps}).

The fundamental group of the Whitehead link exterior admits the presentation\vspace*{-2ex}
\begin{equation}\label{eq:gr-pres-W}
  \pi_{1}(W)=
    \big<\  
      \mu_{0},\ \mu_{1}\ |\ 
      \mu_{0}\mu_{1}\mu_{0}^{-1}\mu_{1}^{-1}\mu_{0}^{-1}\mu_{1}\mu_{0}\mu_{1} =
      \mu_{1}\mu_{0}\mu_{1}\mu_{0}^{-1}\mu_{1}^{-1}\mu_{0}^{-1}\mu_{1}\mu_{0}\
    \big>
\end{equation}
which then gives the following presentation for the fundamental group of the
filled manifold:\begin{equation}\label{eq:gr-pres-Wpq}
  \begin{aligned}
    \pi_{1}(\Wpq)=\big<\ 
      \mu_{0},\ \mu_{1}\ |\ 
      \mu_{0}\mu_{1}\mu_{0}^{-1}\mu_{1}^{-1}\mu_{0}^{-1}\mu_{1}\mu_{0}\mu_{1} &=
      \mu_{1}\mu_{0}\mu_{1}\mu_{0}^{-1}\mu_{1}^{-1}\mu_{0}^{-1}\mu_{1}\mu_{0},\\
      \mu_{0}^{p}\lambda_{0}^{q}&=1\ 
    \big>,
    \end{aligned}
\end{equation}
where
\begin{equation}\label{eq:lambda0-word}
  \lambda_{0}=\mu_{0}\mu_{1}\mu_{0}\mu_{1}^{-1}\mu_{0}^{-1}\mu_{1}^{-1}\mu_{0}\mu_{1}\mu_{0}^{-2} = 
  \mu_{1}\mu_{0}\mu_{1}^{-1}\mu_{0}^{-1}\mu_{1}^{-1}\mu_{0}\mu_{1}\mu_{0}^{-1}
\end{equation}
(since it commutes with $\mu_{0}$) and
\begin{equation}\label{eq:lambda1-word}
  \lambda_{1}=\mu_{0}\mu_{1}\mu_{0}^{-1}\mu_{1}^{-1}\mu_{0}^{-1}\mu_{1}\mu_{0}\mu_{1}^{-1}.
\end{equation}

Following \cite[section~5.1.9]{ST2002} (see also \cite[section~5.6]{ST2020} in these proceedings), we first describe a convenient
subset of the representation variety $\RV(W)$: Let $\rho\in \RV(W)$ be a
non-abelian representation. Then $\rho$ can be conjugated to be of the
following normal form:
\begin{equation}\label{eq:normal-form-generators-RVC}
  \rho(\mu_{0})=
  \begin{pmatrix}
    s & c\\
    0 & s^{-1}
  \end{pmatrix}\qquad\text{and}\qquad
    \rho(\mu_{1})=
  \begin{pmatrix}
    u & 0\\
    1 & u^{-1}
  \end{pmatrix}
\end{equation}
for $s,u\in\C^{\ast}$ and $c\in\C$ (note that $c=0$ when $\rho$ is non-abelian
reducible). We will denote the subset of representations in $\RV(\Wpq)$ of
this form by $\RVC(\Wpq)$.

Applying \eqref{eq:normal-form-generators-RVC} to the left-hand side $w_{1}$
and right-hand side $w_{2}$ of the defining relation in \eqref{eq:gr-pres-W}
yields
\begin{equation}
  \label{eq:Wrh-group-relation-irred-slice}
  0=\rho(w_{1})-\rho(w_{2})=
  \begin{pmatrix}
    -cf & c(u-u^{-1})f\\[1ex]
    (s-s^{-1})f & cf
  \end{pmatrix},
\end{equation}
where
\begin{equation}\label{eq:Wrh-defining-equation-irred-slice}
  \begin{gathered}
  f(s,u,c)=
(s-s^{-1})(u-u^{-1}) + c(s^{-2}u^{-2}-u^{-2}-s^{-2}+4-s^2-u^2+s^2u^2) \\
     +c^2(2s^{-1}u^{-1}-su^{-1}-s^{-1}u+2su)+c^3.
  \end{gathered}
\end{equation}
This shows that $\RVC(W)\subset \C^{3}$ is defined by the equation
$f(s,u,c)=0$. 

Let now $\rho\in\RV(\Wpq)$ be a representation which sends the meridian
$\mu_{1}$ to a parabolic element (a matrix $A\neq \pm I$ with $|\trace A|=2$),
and which has non-abelian image in $\SL_{2}(\C)$. Since $\mu_{1}$ and
$\lambda_{1}$ commute, this implies that $\rho(\lambda_{1})$ is either $\pm I$
or also parabolic.  Following Riley (cf.\ \cite{MR0300267} and
\cite{MR0413078}) we call such a representation an $\SL_{2}(\C)$
\emph{p-representation}, or, for short, \emph{p-rep}. (We remark that Riley
  actually used the name p-rep for $\PSL_{2}(\C)$ representations only, whereas a
  $\SL_{2}(\C)$ representation which projects to a $\PSL_{2}(\C)$ p-rep was called
  \emph{sl-rep}, c.f.\ \cite{MR0413078}.)

Denote by $\fP\subset\RVC(\Wpq)$ the subset of p-reps of $\pi_{1}(\Wpq)$ of
the form~\eqref{eq:normal-form-generators-RVC}. Note that by construction,
$\RVC(\Wpq)$ projects generically $4$-to-$1$ to $\CV(\Wpq)$, corresponding to
the action of the Kleinian four-group (generated here by
the inversions $(s,u,c)\mapsto(s^{-1},u^{-1},c)$ and $(s,u,c)\mapsto\left(
  s,u^{-1},c+(s-s^{-1})(u-u^{-1})\right)$) on the set of possible bases for the
normal form \eqref{eq:normal-form-generators-RVC}.   However, for a p-rep
$\rho\in\RVC(\Wpq)$, we have $u^{2}=1$ (and $s^{2}\neq 1$) and hence the
projection is $2$-to-$1$.  Each conjugacy class of p-reps has therefore
precisely two representatives in $\fP$.

Since $\RV(\Wpq)\subset \RV(W)$, each p-rep $\rho\in\RV(\Wpq)$ induces a non-abelian
representation $\rho'\in\RV(W)$ satisfying
$\rho'(\mu_{0}^{p}\lambda_{0}^{q})=I$ which still sends $\mu_{1}$ to a
parabolic element. We will call such a representation $\rho'$ a
\emph{partial} p-rep. Conversely, each partial p-rep $\rho'$ factors
through $\pi_{1}(\Wpq)$, and we therefore have a bijection between the set
$\fP\subset\RVC(\Wpq)$ of p-reps in normal form
\eqref{eq:normal-form-generators-RVC} and the set $\fP'\subset\RVC(W)$ of
partial p-reps in normal form.

Instead of trying to translate the filling relation
$\mu_{0}^{p}\lambda_{0}^{q}=1$ into defining equations for the parameters $s$,
$c$ and $u$ in $\RVC(\Wpq)$ (which is difficult for general $p$ and $q$), we will now
restrict our attention to the eigenvalues of p-reps.


\subsection{Detour through the eigenvalue variety}
\label{sec:detour-through-eigenv}

The eigenvalue variety of a manifold whose boundary
consists of a finite number of tori was introduced by Tillmann in
\cite{MR2104008} as a generalization of the $A$-polynomial for knots, and an
explicit parametrization of the eigenvalue variety of the Whitehead link
exterior can be found in his thesis (cf.\ \cite[section 5.4.3]{ST2002}) and in \cite[Appendix~B.1]{ST2020} in these proceedings.

Let $\rho\in\RVC(W)$. For $\gamma\in\pi_{1}(W)$, let
\begin{equation}
  \label{eq:upper-left-entry-map-Wrh}
  \begin{aligned}
    \theta_{\gamma}\map &\RVC(W)\ra\C\\ 
    &\rho \mapsto \rho(\gamma)_{1,1}
  \end{aligned}
\end{equation}
be the map taking $\rho$ to the upper left entry of the image $\rho(\gamma)$.
Now define the eigenvalue map  to be
\begin{equation}
  \label{eq:eigenvalue-map-Wrh}
  \begin{aligned}
    e\map &\RVC(W)\ra \C^{4}\\ 
       &\rho\mapsto \left( \theta_{\mu_{0}}(\rho),  \theta_{\lambda_{0}}(\rho),
         \theta_{\mu_{1}}(\rho),  \theta_{\lambda_{1}}(\rho) \right). 
  \end{aligned}
\end{equation}
Note that this is well-defined, as each $\rho\in\RVC(W)$ is triangular on
the peripheral subgroup.

The image $e(\RVC(W))\subset \C^{4}$ is a subset of the eigenvalue
variety which we will denote by $\EV_{0}(W)$.
 Tillmann showed in \cite{ST2002} that the map $e$ has degree one, and that furthermore:
\begin{Lemma}[Lemma 5.6 in \cite{ST2002} or Lemma~17 in \cite{ST2020}]
  The varieties $\RVC(W)$ and $\EV_{0}(W)$ are birationally equivalent.
\qed
\end{Lemma}
In particular, the map $e$ has the following inverse map for $s\neq \pm 1$ and $u^{2}\neq s^{2}$:
\begin{align}\label{eq:inverse-eigenvalue-map}
  e^{-1}&: \EV_{0}(W)\ra \RVC(W), \qquad (s,t,u,v)\mapsto \left(s,u,\frac{s^{2}(t-1)+u^{2}(1-v)}{s^{-1}u^{-1}(s^{2}-u^{2})}\right)\\
  \intertext{As shown in \cite{ST2002}, one can in fact find a pair of inverse maps
  which together are defined on all
  but eight points of $\EV_{0}(W)$:}
  \label{eq:inverse-pair-of-eigenvalue-maps}
  e_{1}^{-1}&: \EV_{0}(W)\ra \RVC(W), \qquad (s,t,u,v)\mapsto \left(s,u,\frac{(s^{2}-v)(1-u^{2})}{su(1+v)}\right)\\
  e_{2}^{-1}&: \EV_{0}(W)\ra \RVC(W), \qquad (s,t,u,v)\mapsto \left(s,u,\frac{(u^{2}-t)(1-s^{2})}{su(1+t)}\right)
\end{align}
Four of the points which are not defined by these two maps correspond to the
discrete faithful representations of $\pi_{1}(W)$ associated to the complete
hyperbolic structure of $W$, parametrized by $t=v=-1$ and
$s^{2}=u^{2}=1$. The remaining four points are given by $t=v=-1$ and
$s^{2}=u^{2}=-1$.

Since a partial p-rep $\rho'\in \RVC(W)$ satisfies $u^{2}=1$ and
$s^{2}\neq 1$, we see that the eigenvalue map $e$ sends the set $\fP'$
bijectively onto a set $\fP''\subset \EV_{0}$, as shown in the following diagram:
\begin{equation}
\label{eq:diagram-slice-and-eigenvalue-varieties}
\xymatrix{
 \RVC(\Wpq) \ar[r]^{\subset} & \RVC(W)    \ar[r]^{\biratequiv} & \EV_{0}(W)\\
 \fP \ar[u] \ar[r]^{\isom}   & \fP'\ar[u] \ar[r]^{\isom}       & \fP'' \ar[u]
}
\end{equation}
We therefore have a bijective correspondence between the set $\fP$ of p-reps in
normal form and a set $\fP''$ of eigenvalue tuples in $\EV_{0}$.
Since each p-rep $\rho\in\RVC(\Wpq)$ and each partial p-rep $\rho'\in\RVC(W)$ satisfies
$\rho(\mu_{0}^{p}\lambda_{0}^{q})=I$, we obtain that the 
set $\fP''$ lies in the intersection $\fF$ of
the following subvarieties of $\EV_{0}(W)$:%
\begin{equation}
  \fP''\subsetneq \fF:=\EV_{0}(W)\intersect \{s^{p}t^{q}=1\} \intersect
  \{u^{2}=1\}.
\end{equation}
It is important to note that $\fP''$ is a proper subset of $\fF$, as
there are indeed points in $\EV_{0}(W)\intersect \{s^{p}t^{q}=1\}$
corresponding to representations which do not factor through $\pi_{1}(\Wpq)$
(we will see that these are precisely points which satisfy $s^{2}=1$).

In summary, we have reduced the generally hard problem of determining the set
$\fP$ to the problem of determining the set $\fP''$. This task is made feasible by
the fact that we can first intersect $\EV_{0}(W)$ with the set $\{u^{2}=1\}$
(which will significantly simplify the defining polynomials involved), and
\emph{then} intersect the resulting set with the subset $\{s^{p}t^{q}=1\}$
corresponding to the Dehn filling.

\subsection{The subset of the eigenvalue variety corresponding to p-reps}
\label{sec:subv-eigenv-variety}

Using the defining polynomial \eqref{eq:Wrh-defining-equation-irred-slice} of
$\RVC(W)$ and using elementary elimination theory, one finds that
$e(\RVC(W))=\EV_{0}(W)\subset \C^{4}$ can be defined as an algebraic set by
the following three polynomials (cf.\ \cite[section 5.4.3]{ST2002} or \cite[Appendix B.1]{ST2020}):
\begin{subequations}
  \label{eq:param-E0-h123}
  \begin{align}
    h_1 &= 
    \begin{aligned}[t]
      &t - s^2 t+s^2 t^2-s^4 t^2-u^2-2 s^2 t u^2+s^4 t u^2-t^2 u^2+2 s^2 t^2 u^2\\
      &+s^4 t^3 u^2+t u^4-s^2 t u^4+s^2 t^2 u^4-s^4 t^2 u^4,
    \end{aligned}\\
    h_2 &=
    \begin{aligned}[t]
      &s^2-v-s^4 v+u^2 v+2 s^2 u^2 v+s^4 u^2 v-s^2 u^4 v+s^2 v^2-u^2 v^2\\
      &-2 s^2 u^2 v^2-s^4 u^2 v^2+u^4 v^2+s^4 u^4 v^2-s^2 u^4 v^3,
    \end{aligned}\\
    h_3 &=
    \begin{aligned}[t]
      &s^4 t-s^6 t-s^2 t u^2+s^4 t u^2+s^6 t^2 u^2-s^2 u^2 v \\
      &+u^4 v+s^2 u^4 v-2 s^4 t u^4 v-u^6 v+s^2 u^6 v^2. 
    \end{aligned}
  \end{align}
\end{subequations}
Intersecting with the subvariety given by $\{u^{2}=1\}$ then simplifies the
defining polynomials to%
\begin{subequations}\label{eq:def-pols-first-intersection}
  \begin{align}
    g_1&= 
    (t-1)\left[ t(t-1)s^{4}+4ts^{2}+(1-t)\right],\\
    g_2&=s^{2}(1-v)(v+1)^{2},\\
    g_3&=s^{2}\left[ t(t-1) s^{4}  + 2t( 1 - v ) s^{2}  + (v^{2}  - t)\right].
  \end{align}
\end{subequations}
Since $s\neq 0$, it can be easily verified that the subvariety defined by these
polynomials splits into the following two irreducible components:%
\begin{equation}
  \label{eq:param-E0-intersect-u2-decomp}
  \begin{gathered}
    \fE_{0}(W)\intersect \{u^{2}=1\} =\\ \{ u^{2}=1,\ t=v=1 \} \union
    \{u^{2}=1,\ v=-1,\ s^{4}t^{2}+(-s^{4}+4s^{2}-1)t+1=0\}
  \end{gathered}
\end{equation}
Note that points in the first component correspond to reducible representations. 

Intersecting this set further, we obtain the set $\fF=\fE_{0}(W)\intersect \{u^{2}=1\}
\intersect \{s^{p}t^{q}=1\}$ as
\begin{equation}
  \label{eq:decomposition-sub-eigenvalue-variety}
  \begin{split}
    \fF = \fF^{r}\union \fF^{i} 
    &= \{u^{2}=t=v=1,\ s^{p}=1\} \\
    & \quad \union 
    \{u^{2}=1, v=-1,\ s^{4}t^{2}+(-s^{4}+4s^{2}-1)t+1=0,\ s^{p}t^{q}=1\}.
  \end{split}
\end{equation}

Let us denote the subsets of $\fP''$ corresponding to reducible and
irreducible p-reps by $(\fP'')^{r}\subsetneq\fF^{r}$ and
$(\fP'')^{i}\subsetneq\fF^{i}$, respectively, so that $\fP''=(\fP'')^{r}\union
(\fP'')^{i}$.


\subsection{Eigenvalues of reducible p-reps}
\label{sec:red-solut}

From equation~\eqref{eq:decomposition-sub-eigenvalue-variety}, we can
immediately read off the points corresponding to reducible
representations as
\begin{equation}
  \label{eq:reducible-solutions}
\fF^{r}=\left\{ (s,1,\pm 1,1) \where s^{p}=1 \right\},
\end{equation}
i.e.\ the eigenvalue $s$ is necessarily a $|p|$-th root of unity.

We now need to determine which of these points also correspond to p-reps
$\rho\in\RV(\Wpq)$, or equivalently to partial p-reps $\rho'\in\RV(W)$.

The reader can easily verify (for example by using the inverse map defined in
\eqref{eq:inverse-pair-of-eigenvalue-maps}) that a reducible representation $\rho'\in\RVC(W)$
corresponding to the eigenvalues $(\pm 1,1,\pm1, 1)$ is necessarily
abelian  and hence
\emph{not} a partial p-rep. On the other hand, all points $(s,1,\pm
1,1)\in\fF^{r}$ with $s\neq\pm 1$ give indeed rise to partial p-reps: Since
$u^{2}=t=v=1$ implies $\rho'(\lambda_{0})=I$, we have
$\rho'(\mu_{0}^{p}\lambda_{0}^{q})=\rho'(\mu_{0}^{p})=I$ and hence $\rho'$
factors through $\pi_{1}(\Wpq)$. Furthermore, $s^{2}\neq 1$ implies that
$\rho'$ is indeed non-abelian.

This means that the subset $(\fP'')^{r}\subsetneq \fF^{r}$ is given by
\begin{equation}
  \label{eq:reducible-solutions-fP''}
  (\fP'')^{r}=\left\{ \left(s,1,\pm 1,1\right) \where s^{p}=1,\ s\neq \pm 1 \right\}.
\end{equation}

To obtain the cardinality of this set, first note that $-1$ is a
$|p|$-th root of unity iff $p$ is even. Since we always have two choices for
$u=\pm 1$, we therefore obtain
\begin{equation}\label{eq:count-fP''-reducible}
\#(\fP'')^{r}=
\begin{cases}
  2(|p|-1) & \text{if $p$ is odd,}\\
  2(|p|-2) & \text{if $p$ is even.} 
\end{cases}
\end{equation}


\subsection{Eigenvalues of irreducible p-reps}
\label{sec:irred-solut}

Let us denote the defining polynomials of $\fF^{i}$ by
\begin{align}
  k_{1}&=s^{p}t^{q}-1\label{eq:defining-equations-k1},\\
  k_{2}&=s^{4}t^{2}+(-s^{4}+4s^{2}-1)t+1,\label{eq:defining-equations-k2}
\end{align}
so that $\fF^{i}=\{ u^{2}=1,\ v=-1,\ k_{1}=k_{2}=0 \}$.

To determine the solutions of $k_{1}=k_{2}=0$ we will use the resultant of
$k_{1}$ and $k_{2}$ with respect to $t$, which we determine by induction on
$q\geq 1$ as the determinant of the corresponding Sylvester matrix:%
\begin{equation}\label{eq:definition-res-pq}
  \begin{split}
    \res_{p,q}(s):&=\text{Resultant}(k_{1},k_{2},t)\\
    &=\det{\footnotesize \begin{pmatrix}
      s^{p}  & 0     & s^{4}           &                 &      & \\
      0      & s^{p} & -s^{4}+4s^{2}-1 & s^{4}           &      & \\
      & 0     & 1               & -s^{4}+4s^{2}-1 &      & \\
      &       &                 & 1               &\ddots& \\
      \vdots & \vdots&                 &                 &\ddots& \\
      &       &                 &                 &\ddots&  s^{4}          \\
      -1     &  0    &                 &                 &      &  -s^{4}+4s^{2}-1\\
      0      &  -1   &                 &                 &      &  1              
    \end{pmatrix}_{(q+2)\times(q+2)}}\\
    &=s^{p+2q} \left( s^{p-2q} + (-1)^{q+1}2T_{q}(y) + s^{-p+2q}
    \right)\\
    &\doteq s^{p-2q} + (-1)^{q+1}2T_{q}(y) + s^{-p+2q},
  \end{split}
\end{equation}
where $y=y(s)=\tfrac{1}{2}(-s^{2}+4-s^{-2})$, and $T_{q}(y)$ is the $q$-th
Chebyshev polynomial of first type, satisfying $T_{0}(y)=1$, $T_{1}(y)=y$ and
$T_{k+1}(y)=2y\ T_{k}(y)-T_{k-1}(y)$ for $k\geq 1$ (see standard literature on
orthogonal polynomials, e.g.\ \cite{MR0372517}).

Note that $\res_{p,q}$ is in fact a Laurent polynomial in $\Z[s^{\pm 1}]$ when
$p<0$, but we will often simply use the word polynomial to also
include Laurent polynomials. Since we are only interested in solutions with
$s\neq 0$, we will always regard this polynomial as being defined up to
multiplication by a unit $s^{k}$ ($k\in \Z$), which we indicate with the
symbol $\doteq$. Since $T_{q}(y)$ is a polynomial of degree $q$ in $y$, the
middle terms of $\res_{p,q}$ are formed by a symmetrical Laurent polynomial of
span $4q$ in the variable $s$, containing only even powers. Note that we can
consider equation~\eqref{eq:definition-res-pq} as a formal definition for the
case $\pqslope=1/0$, as then $\res_{1,0}\doteq (s-1)^{2}=0$ is equivalent to
$k_{1}=0$, implying $s=1$, $t=-1$.

For us, the main importance of this resultant is that its roots characterize
the solutions of the system $\{k_{1}=k_{2}=0\}$: The extension theorem of
elimination theory (cf.\ \cite[\S 1, Theorem 3]{MR1417938}) shows that the
system $\{k_{1}=0,\ k_{2}=0\}$ has a solution $(s,t)$ if and only if $s\neq 0$
is a root of this resultant. We will call $s$ a \emph{partial} solution in
this case.

Furthermore note that given a partial solution
$s$, the value for $t$ is \emph{uniquely} determined by $k_{1}=0$ and
$k_{2}=0$: To see this, note that $k_{2}=0$ implies $t^{2}=-s^{-4}\left(
  (-s^{4}+4s^{2}-1)t+1 \right)$, and hence higher powers of $t$ can
recursively be replaced by linear terms in $t$. The equation
$k_{1}=s^{p}t^{q}-1=0$ can therefore always be expanded to an equation which is
linear in $t$.

Denote the unique solution for $t$ for a given $s$ by $t_{s}$ --- each root
$s\neq 0$ of the resultant therefore corresponds to two eigenvalue points $(s,t_{s},\pm
1,1)$ in $\fF^{i}\subset\fE_{0}(W)$, so that we have
\begin{equation}
  \label{eq:summary-fF-irred}
  \fF^{i}=\{(s,t_{s},\pm 1, 1) \where \res_{p,q}(s)=0 \}.
\end{equation}

As in the reducible case, not all of these points correspond to
representations which will factor through $\pi_{1}(\Wpq)$ to give p-reps in
$\RV(\Wpq)$. This is due to the fact that the eigenvalue condition
$s^{p}t^{q}=1$ does not necessarily imply that
$\rho(\mu_{0})^{p}\rho(\lambda_{0})^{q}=I$. In fact, we have:%
\begin{Proposition}
  The roots $s\notin \{0,1,-1\}$ of $\res_{p,q}$ correspond to irreducible
  p-reps of $\pi_{1}(\Wpq)$; the roots $s=\pm 1$ (when they arise) correspond
  to discrete faithful representations of $\pi_{1}(W)$ which therefore do
  not factor through $\pi_{1}(\Wpq)$.
\end{Proposition}
\begin{proof}
  For $s=\pm 1$, a direct calculation
  shows that an irreducible representation with $s^{2}=u^{2}=1$ is
  conjugate to the form
  \begin{equation}
    \rho(\mu_{0})=\begin{pmatrix}s&-su\pm i\\0&s\end{pmatrix},\ 
    \rho(\mu_{1})=\begin{pmatrix}u&0\\1&u\end{pmatrix},
    \quad\text{and}\ 
    \rho(\lambda_{0})=\begin{pmatrix}-1&-4u\\0&-1\end{pmatrix},
  \end{equation}
  where $s,u\in\{\pm 1\}$. (These are in fact, up to conjugation, all 
  discrete faithful representations into $\SL_{2}(\C)$ corresponding to the
  unique complete hyperbolic structure of the unfilled Whitehead link
  complement.) For each choice of $s=\pm 1$, there are four such representations.
  As the subgroup generated by the parabolic elements
  $\rho(\mu_{0}),\rho(\lambda_{0})\in\SL_{2}(\C)$ is therefore isomorphic to
  $\Z^{2}$, the equation $\rho(\mu_{0})^{p}\rho(\lambda_{0})^{q}=I$ can only
  be satisfied by $p=q=0$ (which is not a valid filling slope).

  For $s\neq \pm 1$, first note that the inverse map
  \eqref{eq:inverse-eigenvalue-map} of the eigenvalue map is well defined
  (since $u^{2}=1$). The two points $(s,t_{s},\pm 1,-1)$ associated to a
  root $s$ therefore correspond to two partial p-reps in $\RVC(W)$.  To see
  that these indeed factor through $\pi_{1}(\Wpq)$, note that for $s\neq \pm
  1$, both $\rho(\mu_{0})$ and $\rho(\lambda_{0})$ are simultaneously
  diagonizable, since they are not parabolic and commute. The condition
  $s^{p}t^{q}=1$ then immediately implies
  $\rho(\mu_{0})^{p}\rho(\lambda_{0})^{q}=I$.
\end{proof}

This concludes the proof of
Theorem~\ref{thm:conjugacy-classes-preps-charact-by-roots-of-respq}. For the
following application to the minimal total Culler-Shalen seminorm in
section~\ref{sec:minimal-CS-norm}, we now need to study the number of conjugacy
classes of p-reps.


\subsection{Roots of the resultant}
\label{sec:roots-resultant}

We summarize some key facts about the roots of the resultant $\res_{p,q}$:

\begin{Lemma}[Trivial roots]\label{lem:properties-roots-respq}
  If $q$ is even (and hence $p$ odd), $s=1$ is a root of order $2$ of
  $\res_{p,q}$. If both $q$ and $p$ are odd, $s=-1$ is a root of order $2$ of
  $\res_{p,q}$. In the remaining case ($q$ odd and $p$ even), neither of $\pm
  1$ is a root of $\res_{p,q}$.
\end{Lemma}
\begin{proof}
  This can be directly verified from the equations $k_{1}=k_{2}=0$ (cf.\
  \eqref{eq:defining-equations-k1} and~\eqref{eq:defining-equations-k2}): If
  $s=\pm 1$, equation $k_{2}=0$ necessarily implies $t_{s}=-1$.  Using
  $k_{1}=0$, we then immediately see that $(1,-1)$ is a solution iff $q$ is
  even, and that $(-1,-1)$ is a solution iff $p+q$ is even (which --- under
  the assumption that $p$ and $q$ are coprime --- is equivalent to both $p$
  and $q$ being odd). To show that the order is equal to $2$, note that
  \begin{align*}
    \res'_{p,q}&:= \frac{d}{ds} \res_{p,q}= 
    (p-2q)\left( s^{p-2q-1}-s^{-p+2q-1}\right) + (-1)^{q+1}2q U_{q-1}(y) y'\\
    \intertext{and}
    \res''_{p,q}&= 
    \begin{gathered}[t]
      (p-2q)\left( (p-2q-1)s^{p-2q-2}+(-p+2q-1)s^{-p+2q-2}\right)\\
      + (-1)^{q+1}2q \left( y'' U_{q-1}(y) + (y')^{2} \tfrac{\der}{\der y}U_{q-1}(y) \right),
    \end{gathered}%
  \end{align*}
  where $U(y)$ is the Chebyshev polynomial of second type, satisfying
  $U_{0}(y)=1$, $U_{1}(y)=2y$ and $U_{k+1}(y)=2y U_{k}(y)-U_{k-1}(y)$ for
  $k\geq 1$ (cf.\ \cite{MR0372517}). Since $y'(\pm 1)=0$, we immediately get
  $\res'_{p,q}(\pm 1)=0$. Furthermore, $\res''_{p,q}(\pm 1)\neq 0$ unless
  $p/q=0/1$ or $p/q=4/1$, but in those cases $\res_{p,q}$ degenerates and has
  no non-zero roots.
\end{proof}

\begin{Lemma}[Symmetries]~\\[-2.5ex]\label{lem:symmetries-res}
  \begin{enumerate}
  \item $\res_{p,q}(s)=0$ $\iff$ $\res_{p,q}(s^{-1})=0$.\label{lem:sym-res-inv}
  \item $\res_{p,q}(s)=0$ $\iff$ $\res_{p,q}(\bar{s})=0$.\label{lem:sym-res-conj}
  \item $(\res_{p,q}(s)=0$ $\iff$ $\res_{p,q}(-s)=0)$ iff $p$ is even.\label{lem:sym-res-neg}
  \item $\res_{p,q}(s)\doteq \res_{-p+4q,q}(s)$.\label{lem:sym-res-negate-p}
  \end{enumerate}
\end{Lemma}
\begin{proof}
  For (1), observe that $y(s^{-1})=y(s)$. For (2), note that $\res_{p,q}$ is a
  Laurent polynomial with real (in fact integer) coefficients. For (3), first
  note that $y(-s)=y(s)$ and hence \[\res_{p,q}(-s)-\res_{p,q}(s)=
  \left((-1)^{p}-1\right) \left(s^{p-2q}+s^{-p+2q}\right).\] This is clearly
  $0$ when $p$ is even. Conversely, assume this difference to be $0$ but $p$
  to be odd. Then $s^{2(p-2q)}=-1$, which implies $s=e^{i\phi}$ for some
  $0\leq \phi<2\pi$. But then
  $\res_{p,q}(s)=(-1)^{q+1}2T_{q}(y(s))=(-1)^{q+1}2T_{q}(2-\cos(2\phi))\neq
  0$, since the roots of $T_{q}(y)$ lie in the interval $(0,1)$ (cf.\
  \cite{MR0372517}).
  The last statement is obvious from \eqref{eq:definition-res-pq}.
\end{proof}

\begin{Lemma}[Real roots]\label{thm:res-real-roots}
  When $p>4q>0$ or $p<0$, the only possible real-valued roots of $\res_{p,q}$ are $\pm
  1$.  When $0<p<4q$ and $p$ is odd, there are two additional real-valued, positive
  roots. When $0<p<4q$ and $p$ is even (and hence $q$ odd), there
  are precisely four real-valued roots in total, none of which is $\pm 1$.
\end{Lemma}
\begin{proof}
  Assume $\res_{p,q}(s)=0$ with real-valued $s>0$. Then
  \[s^{p-2q}+s^{-(p-2q)}=(-1)^{q}2T_{q}(y(s)).\]  Denote the left-hand side of
  this equation by $\phi(s)$, the right-hand side by $\psi(s)$. We will study
  intersections of these two real-valued functions. Note that $\phi(s)$ is
  strictly monotonously decreasing for $0<s<1$, has a global minimum of value
  $2$ at $s=1$, and is strictly monotonously increasing for $s>1$.
  Furthermore, note that the range of $y(s)$ is $(-\infty,1]$, with global
  maximum at $s=1$. This implies that the range of
  $\psi(s)=(-1)^{q}2T_{q}(y(s))$ is $(-\infty,2)$, using the fact that $T_{q}(y)$ is a
  polynomial of degree $q$ in $y$, whose zeros and relative extrema of value
  $\pm 1$ are all contained in the interval $[-1,1]$. More specifically,
  $\psi$ oscillates with values between $-2$ and $2$ for $s\in
  [1-\sqrt{2},1+\sqrt{2}]$ (with $\psi(1)=(-1)^{q}2$) and tends
  strictly monotonously to $+\infty$ for $s\ra 0$ and $s\ra \infty$.
  
  Apart from the possibility of an intersection at $s=1$ (in the case that $q$
  is even), other intersections between $\phi$ and $\psi$ can therefore only
  occur outside the interval $(1-\sqrt{2},1+\sqrt{2})$ when the maximal degree
  in $s$ of $\psi$ is larger than that of $\phi$.  But $\maxdeg_{s}
  \phi=|p-2q|=-\mindeg_{s} \phi$ and $\maxdeg_{s} \psi=2q=-\mindeg_{s} \psi$.
  It is then easy to see that for $p>4q>0$ or $p<0$, $\maxdeg_{s}
  \psi<\maxdeg_{s} \phi$ and there are no such intersections. The same
  argument then shows that for $0<p<4q$, there are precisely two positive
  real-valued roots different from $1$ (which are inverses of each other, by
  Lemma~\ref{lem:symmetries-res}~(\ref{lem:sym-res-inv})).
  
  To analyze the remaining case $s<0$, note that $\phi(-s)=(-1)^{p}\phi(s)$
  while $\psi(-s)=\psi(s)$. If $p$ is odd, this means that there cannot be any
  intersection apart from possibly $s=-1$ (which arises if $q$ is odd), as the
  two functions diverge in opposite directions. If $p$ is even (and hence $q$
  odd), we get as roots the negatives of the roots found for $s>0$. In total,
  there are either $0$ or $4$ real roots in this case, as $\res_{p,q}(\pm
  1)\neq 0$.
\end{proof}

Note that this includes the case $p/q=2/1$, when there are only four non-zero roots in
total, all of which are real and distinct from $\pm 1$.

\begin{Lemma}[Imaginary roots]\label{thm:res-imaginary-roots}
  When $p>4q>0$ or $p<0$,  $\res_{p,q}$ has imaginary roots iff $p\in 4\Z$. In
  that case, there are four such roots.
  When $0<p<4q$, there are no imaginary roots.
\end{Lemma}
\begin{proof}
  Similar to the proof of Lemma~\ref{thm:res-real-roots}. Assume 
  $\res_{p,q}(ix)=0$ with $x>0$. Then
  \[i^{-p+2q}\left((-1)^{p}x^{p-2q}+x^{-p+2q}\right)
  =(-1)^{q}2T_{q}\left(\tfrac{x^{2}+4+x^{-2}}{2}\right).\]
  Now note that
  the right-hand side is always real-valued, while the same for the left-hand
  side is only true when $p$ is even. In that case, the equation
  is equivalent to
  \[\phi(x):=(-1)^{p/2}\left(x^{p-2q}+x^{-p+2q}\right)
  =2T_{q}\left(\tfrac{x^{2}+4+x^{-2}}{2}\right)=:\psi(x).\]
  As $\psi(x)>0$ for all $x\in \R$, this equation can only be satisfied if $p/2$
  is even, i.e.\ if $p\in 4\Z$. Under this assumption, comparison of the
  maximal degrees then shows that there is no intersection of $\phi$ and
  $\psi$ when $0<p<4q$, and precisely two intersections (which are inverses of
  each other) when $p>4q>0$ or $p<0$.  The symmetry with respect to complex
  conjugation then shows that there are in total four pure imaginary roots in this
  case, and none when $p\notin 4\Z$.
\end{proof}

\begin{Lemma}\label{lem:res-only-roots-on-unit-circle-are-pm1}
 None of the roots $s\notin\{\pm 1\}$ of $\res_{p,q}$ lies on the unit circle.
\end{Lemma}
\begin{proof}
  Simply note that for $\phi\in[0,2\pi)$,
\begin{equation*}
  \res_{p,q}(e^{i\phi}) =2(\underbrace{\cos(p-2q)\phi}_{\in[-1,1]} + (-1)^{q+1}\underbrace{T_{q}(2-\cos
  2\phi)}_{\geq 1})=0
\end{equation*}
implies $\phi=0$ or $\phi=\pi$.
\end{proof}

Examples of the typical distribution of the roots are shown in
figures~\ref{fig:res-roots-typical-distribution1},
\ref{fig:res-roots-typical-distribution2} and
\ref{fig:roots-of-res-examples-close-to-4} on page~\pageref{fig:res-roots-typical-distribution1}.

\begin{figure}[htb]
  \centering
  \begin{minipage}[b]{0.45\textwidth}%
    \centering
    \subfigure[$\res_{65,3}$ ($p>4q>0$)]{%
      \centering\includegraphics[width=0.6\textwidth]{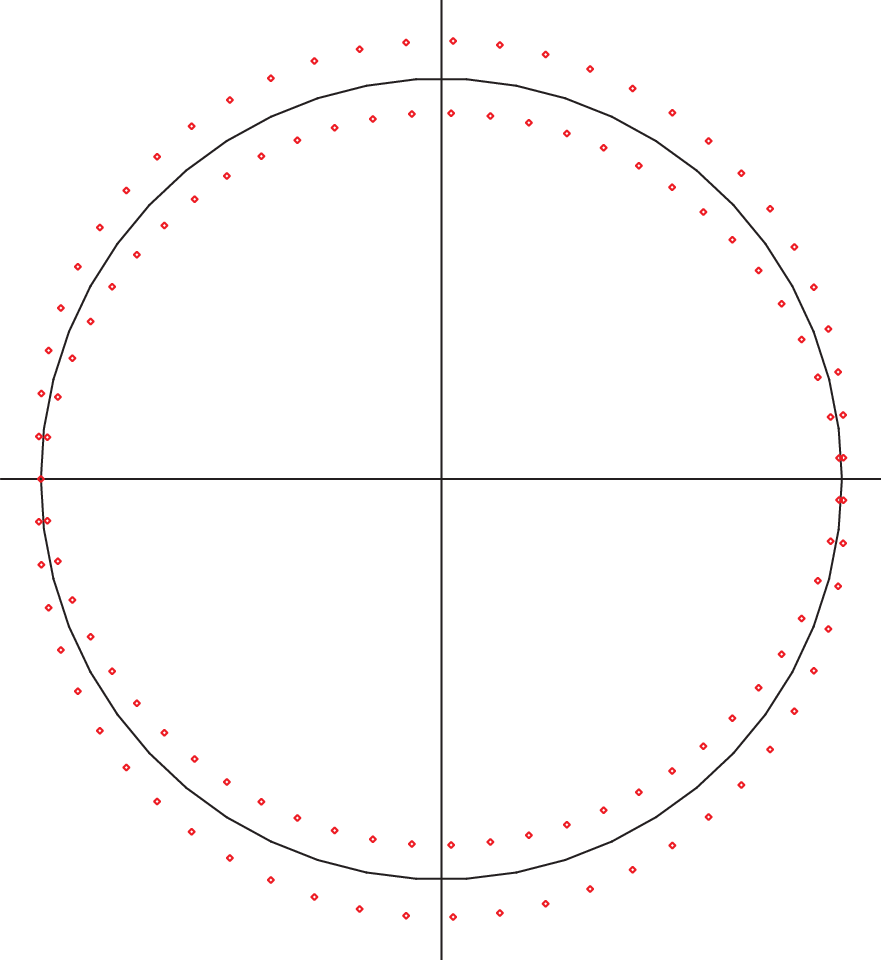}%
      \label{fig:res-roots-typical-distribution1}%
    }\vspace{1cm}
    \subfigure[$\res_{65,23}$ ($0<p<4q$)]{%
        \centering\includegraphics[width=\textwidth]{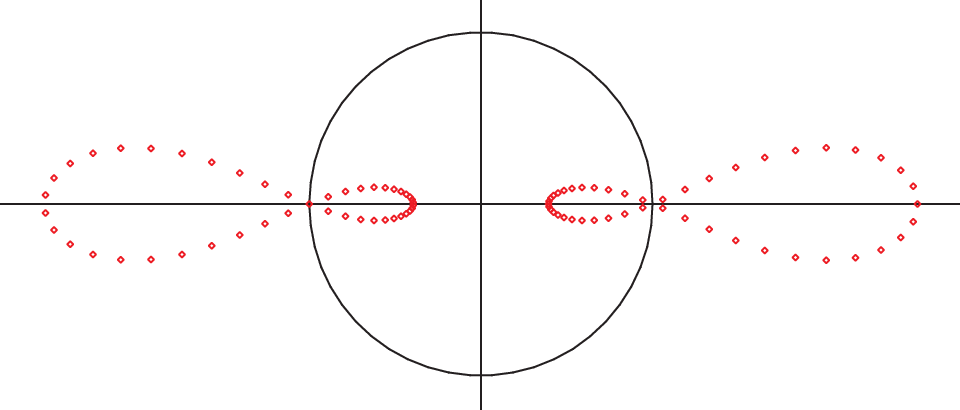}%
        \label{fig:res-roots-typical-distribution2}%
    }
  \end{minipage}
  \begin{minipage}[b]{0.45\textwidth}%
    \subfigure[$\res_{65,16}$ ($p-4q>0$, but small)]{%
      \includegraphics[width=\textwidth]{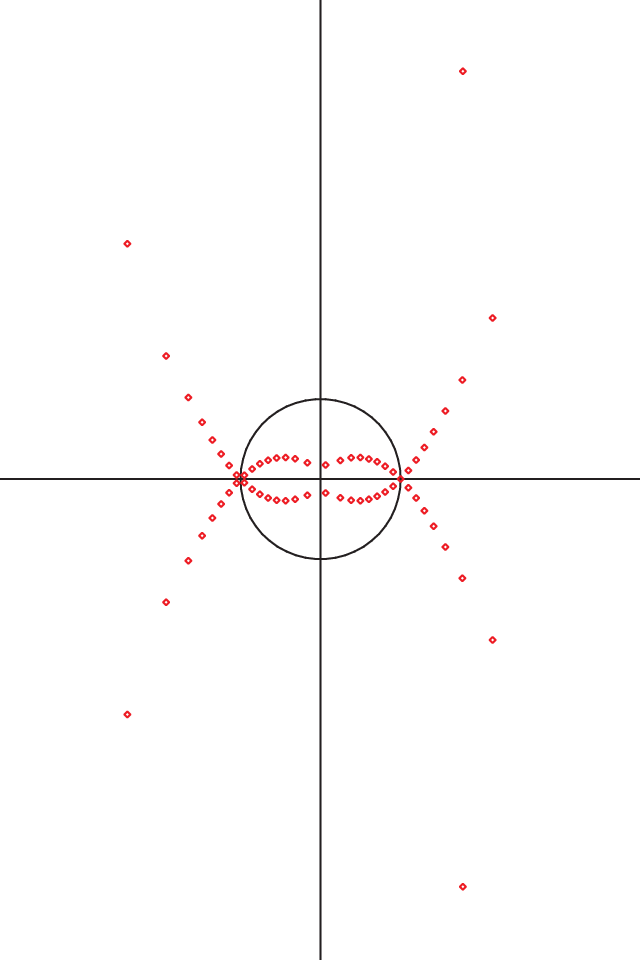}%
      \label{fig:roots-of-res-examples-close-to-4}%
    }
  \end{minipage}
  \caption{Examples of the distribution of the non-zero roots of $\res_{p,q}$}
  \label{fig:roots-of-res-examples-generic}
\end{figure}

We will now determine the number of non-zero roots of the resultant.  The
number of non-zero roots (with multiplicities) of a Laurent polynomial is
given by its span, i.e.\ the difference of the maximal and minimal degrees.
From \eqref{eq:definition-res-pq}, we see that unless $\pqslope\in \{0,4\}$,
the span of $\res_{p,q}$ is $2\max(|p-2q|,2q)$ (recall that we always assume
$q>0$). When $\pqslope\in \{0,4\}$, $\res_{p,q}\doteq 1$, and there are
therefore no irreducible p-preps --- matching the fact
that $W_{0}$ and $W_{4}$ are toroidal manifolds which do not admit a complete
hyperbolic structure.

Discarding $0$, $1$ and $-1$ from the roots of $\res_{p,q}$, we therefore
obtain an upper bound for the number of \emph{distinct} non-trivial roots as given in
table~\ref{tab:number-of-simple-roots-resultant}. We remark that it is a hard
problem to analytically show that all non-trivial roots are in fact simple ---
but the proof of Theorem~\ref{thm:total-cs-norm-Wpq} will in fact imply 
this for odd $p$ (cf.\ Corollary~\ref{cor:roots-of-res-are-simple}).
\begin{table}[htbp]
  \centering
  \begin{tabular}{c||CCC}
  $\pqslope\in$ &  (-\infty,0)  & (0,4) & (4,\infty)\\\hline
  $p$ odd &    2|p|+4|q|-2  &   4|q|-2  &  2|p|-4|q|-2 \\
  $p$ even &   2|p|+4|q|    &   4|q|    &  2|p|-4|q| 
  \end{tabular}
  \\[2ex]\caption{Upper bound for the number of distinct roots $s\not\in \{0,\pm 1\}$ of $\res_{p,q}$}
  \label{tab:number-of-simple-roots-resultant}
\end{table}

\subsection{The number of conjugacy classes of p-reps}
\label{sec:number-p-reps}

We will now put together the results of the previous sections to obtain the
count of conjugacy classes of p-reps for the manifold $\Wpq$.

We have determined the two components of the set $\fP''$ of eigenvalue tuples
corresponding to p-reps as
\begin{align}
    (\fP'')^{r} &= \{(s,\pm 1, 1, 1) \where s^{p}=1,\ s\neq \pm 1 \}, \\
    (\fP'')^{i} &= \{(s,\pm 1, t_{s}, -1) \where \res_{p,q}(s)=0,\ s\neq \pm 1 \},
\end{align}
with their cardinalities given in equation~\eqref{eq:count-fP''-reducible}
and table~\ref{tab:number-of-simple-roots-resultant}, respectively.

Now recall that by construction, the
cardinality of $\fP''$ is equal to the cardinality of $\fP'$ and that of
$\fP$. Furthermore, we saw that each conjugacy class of p-reps in $\RVC(\Wpq)$
is represented by \emph{two} elements in $\RVC(\Wpq)$ (corresponding to the
pair of eigenvalue points $(s,t,u,v)$ and $(s^{-1},t^{-1},u,v)$). To obtain
the number of conjugacy classes of p-reps, we therefore have to divide the
cardinality of $\fP$ by $2$.
Since the two representatives for a conjugacy class in $\RVC(\Wpq)$ clearly
have the same character independent of whether they are reducible or
irreducible, we see that this number also coincides with the \emph{number of
  characters} of p-reps.

This then shows:
\begin{Proposition}\label{thm:upper-bounds-for-conjugacy-classes-of-preps}
  The number of distinct conjugacy classes of $\SL_{2}(\C)$ p-reps of
  $\pi_{1}(\Wpq)$ is bounded above by the numbers shown in
  table~\ref{tab:number-of-conjug-classes-preps-Wpq}:
\end{Proposition}

\begin{table}[htbp]
  \centering
  \begin{tabular}{c||CCC}
  $\pqslope\in$ &  (-\infty,0)  & (0,4) & (4,\infty)\\\hline
  $p$ odd &    3|p|+4|q|-3  &   |p|+4|q|-3  &  3|p|-4|q|-3 \\
  $p$ even &   3|p|+4|q|-2  &   |p|+4|q|-2  &  3|p|-4|q|-2 
  \end{tabular}
  \\[2ex]\caption{Upper bound for the number of conjugacy classes of $\SL_{2}(\C)$ p-reps of $\pi_{1}(\Wpq)$}
  \label{tab:number-of-conjug-classes-preps-Wpq}
\end{table}



\newpage
\section{The minimal total Culler-Shalen seminorm and p-reps}
\label{sec:minimal-CS-norm}

The goal of this section is to show:
\begin{importantTheorem}\label{thm:minimal-seminorm-is-number-of-characters-of-preps}
  The minimal total Culler-Shalen seminorm of the manifold $\Wpq$ is equal to
  the number of characters of p-reps of $\pi_{1}(\Wpq)$. In particular, it is
  less than or equal to the bounds given in
  Table~\ref{tab:number-of-conjug-classes-preps-Wpq}. 
\end{importantTheorem}

Since Dehn filling along the meridian $\mu_{1}$ of $\Wpq$ gives
a closed manifold with cyclic fundamental group, while $\mu_{1}$ ($=\frac{1}{0}=\infty$) is not a boundary
slope, Corollary 1.1.4 of \cite{CGLS} implies that its Culler-Shalen seminorm is
in fact equal to the minimal total Culler-Shalen seminorm $s$:
\begin{equation}
  \label{eq:CS-of-mu1-is-s}
  s=\CSnorm{\mu_{1}}=\CSnorm{\infty}.
\end{equation}
We can therefore find $s$ by explicitly determining $\CSnorm{\mu_{1}}$.

Recall that $\CSnorm{\mu_{1}}=\deg (f_{\mu_{1}}\map
\CVsp^{\text{irr}}(\Wpq)\ra \CP^{1})$. We will determine this degree as a
count of the roots of $f_{\mu_{1}}$. First note that we can in fact
restrict our attention to $\CV(\Wpq)$, since $f_{\mu_{1}}$ has poles at each
ideal point in $\CVsp(\Wpq)$ (see e.g.\ the remarks after Corollary 9.2.7 in
\cite{MR1886685}).

Since $f_{\mu_{1}}=0$ is equivalent to $\chi_{\rho}(\mu_{1})=\trace
\rho(\mu_{1})=\pm 2$, we immediately see that zeroes of this function are
characters of representations $\rho\in\RV(\Wpq)$ which send the meridian
$\mu_{1}\in\pi_{1}(\Wpq)$ to $\pm I$ or a parabolic element of
$\SL_{2}(\C)$. Since we only regard non-trivial curves in $\CV(\Wpq)$, the
representations in question have to be non-abelian, and are therefore \emph{p-reps}.

Theorem~\ref{thm:minimal-seminorm-is-number-of-characters-of-preps} will
therefore be a direct consequence of the following two propositions, which are
the subject of the bulk of this section:%

\begin{Proposition}\label{thm:char-of-prep-is-smooth-point}
  Each character of a p-rep of $\pi_{1}(\Wpq)$ is a smooth point of a unique
  non-trivial curve in $\CVi(\Wpq)$.   
\end{Proposition}
\begin{Proposition}\label{thm:char-of-prep-is-simple-zero}
  Each character of an irreducible p-rep of $\pi_{1}(\Wpq)$ is a simple zero
  of $f_{\mu_{1}}$. If $p/q\neq 0$, also each character of a non-abelian
  reducible p-rep of $\pi_{q}(\Wpq)$ is a simple zero of $f_{\mu_{1}}$.
\end{Proposition}

We will begin in section~\ref{sec:smooth-char-non-abelian-reducible-prep} by
using work of Heusener and Porti (cf.\ \cite{MR2171800}) to show
Proposition~\ref{thm:char-of-prep-is-smooth-point} for a \emph{general
  non-abelian reducible} representation (therefore including the case of
reducible p-reps).  We remark that this is a subtle question, as it asks when
the character of a non-abelian reducible representation can be deformed into
the character of an irreducible representation.  This result will also be
needed in section~\ref{sec:whitehead-seifert-fillings}.

To show Proposition~\ref{thm:char-of-prep-is-smooth-point} for an
\emph{irreducible} p-rep, we then exhibit a useful partially-diagonal normal
form for p-reps in section~\ref{sec:param-irred-repr}, and use this in
section~\ref{sec:cohomology} to determine an explicit presentation of the
cohomology of $\pi_{1}(\Wpq)$ with coefficients in $\sltwist$ (the Lie algebra
$\sll_{2}(\C)$ viewed as a $\pi_{1}(\Wpq)$ module via the map $\Ad\after\rho$
for a given p-rep $\rho$).  The knowledge of the dimension of this cohomology
group then implies our claim using Theorem A of \cite{MR1909000}.

Finally, Proposition~\ref{thm:char-of-prep-is-simple-zero} will be shown in
section~\ref{sec:deform-repr-degr} by studying the degree of the trace
function $f_{\mu_{1}}$ in relation to deformations of representations.


\subsection{Smoothness of the characters of general non-abelian reducible representations}
\label{sec:smooth-char-non-abelian-reducible-prep}

We will show:%
\begin{Proposition}\label{thm:Wpq-nonab-red-lie-lie-in-intersection} Let
  $\prho\in\PRV(\Wpq)$ be a non-abelian reducible $\PSL_{2}(\C)$
  representation. For $p/q\neq 4$, the character $\chi_{\prho}$ of $\prho$ is
  contained in precisely two irreducible components of $\PCV(\Wpq)$, both of
  which are curves, one trivial, the other non-trivial. In addition,
  $\chi_{\prho}$ is a smooth point of both curves and the intersection at
  $\chi_{\prho}$ is transverse. When $p/q=4$, the above is true for all
  non-abelian reducible representations $\prho$ with $\trace
  \prho(\mu_{0})\neq \pm 2$.
\end{Proposition}
Since every non-abelian reducible $\SL_{2}(\C)$ representation
$\rho\in\RV(\Wpq)$ induces a non-abelian reducible $\PSL_{2}(\C)$
representation $\prho\in\PRV(\Wpq)$, this then implies the following stronger
version of Proposition~\ref{thm:char-of-prep-is-smooth-point} in the case of
non-abelian reducible representations:%
\begin{Corollary}\label{thm:Wpq-nonab-red-contribute-to-norm}
  Every non-abelian reducible $\SL_{2}(\C)$ character
  $\chi_{\rho}\in\CV(\Wpq)$ is a smooth point on a unique non-trivial curve
  $\CV_{0}\subset\CV(\Wpq)$.
\end{Corollary}
\begin{proof}
  Let $\rho\in\RV(\Wpq)$ be a non-abelian reducible representation, and let
  $\prho\in\PRV(\Wpq)$ be its image under the natural projection $\pi\map
  \RV(\Wpq)\ra\PRV(\Wpq)$.
  Proposition~\ref{thm:Wpq-nonab-red-lie-lie-in-intersection} now shows that
  $\prho$ lies as a smooth point on a unique non-trivial curve
  $\PCV_{0}\subset\PCV(\Wpq)$, and we only need to show that this curve lifts
  to a non-trivial curve $\CV_{0}\subset\CV(\Wpq)$. But this is immediately
  clear, since $\prho$ lifts to $\rho$ by construction, and $\pi$ is a regular
  covering map onto its image (cf.\ \cite[Section~3]{BZseminorms}).
\end{proof}

We start by determining all possible non-abelian reducible representations of
$\pi_{1}(\Wpq)$ up to conjugation:
\begin{Lemma}\label{lem:nonab-red-specific-form-Wpq} Let $\prho$
  be a non-abelian reducible representation in $\PRV(\Wpq)$. Then up to
  conjugation, $\prho$ has the form
  \begin{equation*}
    \prho(\mu_{0})=
    \pm\begin{pmatrix}
      s & 0\\
      c & s^{-1}
    \end{pmatrix}\qquad\text{and}\qquad
    \prho(\mu_{1})=
    \pm\begin{pmatrix}
      u & 0\\
      d & u^{-1}
    \end{pmatrix},
  \end{equation*}
  where either
  \begin{enumerate}
  \item $s^{2}=1$, $c\neq 0$ and $u^{2}=\frac{1}{2q}(-p+2q\pm\sqrt{p(p-4q)})$
    (where $p/q=0$ cannot arise),
    or
  \item $u^{2}=1$, $d\neq 0$ and $s^{p}=\pm 1$, but $s^{2}\neq 1$.
  \end{enumerate}
\end{Lemma}
\begin{proof}
  Let $\prho\in\PRV(\Wpq)$ be non-abelian reducible. Since it is reducible, we
  can conjugate it into lower triangular form as in the claim. The first group
  relation of $\pi_{1}(\Wpq)$ (cf.\ equation \eqref{eq:gr-pres-Wpq}) now
  implies that either $du(1-s^{2})+cs(u^{2}-1)=0$ (in which case $\prho$
  would be abelian, a contradiction), or one of $s^{2}=1$ or $u^{2}=1$ --- but not
  both (otherwise $\prho$ would again be abelian).

  If $s^{2}=1$, we hence have $u^{2}\neq 1$ and $c\neq 0$. We can now conjugate by
  diagonal matrices to let $c=1$. Then $\prho(\lambda_{0})=
  \pm\begin{smallpmatrix}
    1 & 0\\
    s^{-1}u^{-2}(u^{2}-1)^{2} & 1
  \end{smallpmatrix}$, and we get
\[
\pm I=\prho(\mu_{0})^{p}\prho(\lambda_{0})^{q}
\begin{aligned}[t]
&=\pm 
\begin{pmatrix}
  s^{p} & 0\\
  s^{p-1}p & s^{p}
\end{pmatrix}
\begin{pmatrix}
  1 & 0\\
  q s^{-1}u^{-2}(u^{2}-1)^{2} & 1
\end{pmatrix}\\
&=\pm \begin{pmatrix}
  s^{p} & 0\\
  ps^{p-1}+s^{p}q s^{-1}u^{-2}(u^{2}-1)^{2} & s^{p}
\end{pmatrix}.
\end{aligned}
\]
Since $s^{2}=1$, the lower left entry implies $qu^{4}+u^{2}(p-2q)+q=0$, and
therefore $u^{2}=\frac{1}{2q}(-p+2q\pm\sqrt{p(p-4q)})$ as claimed. Note that these
roots are distinct unless $p/q=0$ or $p/q=4$. When $p/q=0$, this implies that
$u^{2}=1$, a contradiction since $\prho$ is non-abelian. When $p/q=4$, we get
$u^{2}=-1$, which yields a valid non-abelian reducible representation.

In the remaining case, when $s^{2}\neq 1$, we have $u^{2}=1$ and can assume
$d=1$. Then $\prho(\lambda_{0})=\pm I$ and we immediately get $\pm
I=\prho(\mu_{0})^{p}\prho(\lambda_{0})^{q}=\prho(\mu_{0})^{p}$ which implies
$s^{p}=\pm 1$.
\end{proof}

Following \cite{MR2171800}, we now determine the twisted Alexander
polynomial associated to the diagonal representation $\prho_{\alpha}$ obtained
from a non-abelian reducible representation $\prho\in\PRV(\Wpq)$ as in
Lemma~\ref{lem:nonab-red-specific-form-Wpq} by setting
all non-diagonal matrix entries to be zero. In particular, this representation has
the same character as $\prho$, i.e.\ $\chi_{\prho_{\alpha}}=\chi_{\prho}$.
\begin{Lemma}\label{lem:tw-alex-nonab-red-Wpq}
  Let $\prho$ be a non-abelian reducible representation in $\PRV(\Wpq)$, and
  let $\prho_{\alpha}$ be the diagonal representation with the same character
  as $\prho$ as described above.  The twisted Alexander
  polynomial associated to $\rho_{\alpha}$ is given by:%
  \begin{equation*}
    \Delta^{\Phi_{\sigma}}(t)\doteq
    \begin{cases}
      t-1 & \text{if $s^{2}\neq 1$,}\\
      qt^{2}+(p-2q)t+q & \text{if $s^{2}= 1$.}
    \end{cases}
  \end{equation*}
\end{Lemma}
\begin{proof}
  We follow the notation of \cite{MR2171800} to facilitate comparison
  for the reader. First note that $H_{1}(\Wpq;\Z)\isom \Zmod{p} \oplus \Z =
  \generate{\mu_{0}} \oplus \generate{\mu_{1}}$, where we use the same symbols
  for the generators of $\pi_{1}(\Wpq)$ and their images under abelianization.
  We now define the following homomorphisms:%
  \begin{equation*}
    \begin{alignedat}[t]{3}
      \alpha& \map\pi_{1}(\Wpq)\ra H_{1}(\Wpq;\Z)\ra \C^{\ast},
            &\qquad& \mu_{0}\mapsto s^{2}, &\quad& \mu_{1}\mapsto u^{2},\\
      \phi& \map\pi_{1}(\Wpq)\ra \Z=\generate{t},
            && \mu_{0}\mapsto 1, && \mu_{1}\mapsto t,\\
      \sigma& \map\pi_{1}(\Wpq)\ra \tors H_{1}(\Wpq;\Z)\overset{\alpha}{\ra} \C^{\ast},
            && \mu_{0}\mapsto s^{2}, && \mu_{1}\mapsto 1,
    \end{alignedat}
  \end{equation*}
  where $s$ and $u$ are the upper left entries of $\prho(\mu_{0})$ and
  $\prho(\mu_{1})$, respectively.
  Using these maps, we can now define a map $\Phi_{\sigma}\map
  \pi_{1}(\Wpq)\ra \C[t^{\pm}]^{\ast}$ to the units of the ring of Laurent
  polynomials $\C[t^{\pm 1}]$ by
  $\Phi_{\sigma}(\gamma):=\sigma(\gamma)t^{\phi(\gamma)}$ for
  $\gamma\in\pi_{1}(\Wpq)$. This map can be linearly extended to the group ring
  $\Z\pi_{1}(\Wpq)$ which allows us to define the twisted Alexander polynomial
  $\Delta^{\Phi_{\sigma}}(t)\in\C[t^{\pm 1}]$ (cf.\ 
  \cite[Section~2]{MR2171800} or e.g.\ \cite{MR1273784} for more
  background). Since we have a
  presentation with two generators $\mu_{0}$ and $\mu_{1}$ and two relators
  \[r_{1}=\mu_{0}\mu_{1}\mu_{0}^{-1}\mu_{1}^{-1}\mu_{0}^{-1}\mu_{1}\mu_{0}\mu_{1}
  \mu_{0}^{-1}\mu_{1}^{-1}\mu_{0}\mu_{1}\mu_{0}\mu_{1}^{-1}\mu_{0}^{-1}\mu_{1}^{-1}\]
  and $r_{2}=\mu_{0}^{p}\lambda_{0}^{q}$, the twisted Alexander polynomial can
  be obtained as the greatest common divisor of the entries of the twisted Fox Jacobian
  $J^{\Phi_{\sigma}}$:%
  \begin{equation*}
    \begin{aligned}
      J^{\Phi_{\sigma}}&=
      \begin{pmatrix}
        \pderx{r_{1}}{\mu_{0}} & \pderx{r_{1}}{\mu_{1}}\\[2ex]
        \pderx{r_{2}}{\mu_{0}} & \pderx{r_{2}}{\mu_{1}}
      \end{pmatrix}^{\!\Phi_{\sigma}}\\
      &=\begin{pmatrix}
        s^{-2}t^{-1}(t-1)^{2}(s^{2}-1) & s^{-2}t^{-1}(1-t)(s^{2}-1)^{2}\\
        \sum_{k=0}^{p-1} s^{2k}+s^{2p}qt^{-1}(t-1)^{2} & s^{2p}q t^{-1}(1-t)(s^{2}-1)
      \end{pmatrix}.
    \end{aligned}
  \end{equation*}
  This matrix can be obtained by hand by first determining the Fox derivatives of the
  group relators (cf.\ \cite[Ch.\ 9B]{MR0053938} or \cite{MR1959408} for
  background on Fox free differential calculus), and then applying the map
  $\Phi_{\sigma}$ to the obtained elements of the group ring
  $\Z\pi_{1}(\Wpq)$.

  We now have two cases: If $s^{2}=1$, we get
  \begin{equation*}
    J^{\Phi_{\sigma}}=
    \begin{pmatrix}
      0 & 0\\
      p+qt^{-1}(t-1)^{2} & 0
    \end{pmatrix},
  \end{equation*}
  and thus that $\Delta^{\Phi_{\sigma}}(t)=p+qt^{-1}(t-1)^{2} \doteq
  qt^{2}+(p-2q)t+q$ (where $\doteq$ denotes equivalence up to multiplication
  by units of $\C[t^{\pm 1}]$). Note that this is in fact the \emph{untwisted}
  Alexander polynomial of $\Wpq$, as $s^{2}=1$ implies that $\sigma$ is the
  trivial $1$-dimensional representation (cf.\ also e.g.\ \cite{MR93i:57019},
  taking into account that the authors work with the \emph{left-handed}
  Whitehead link).

    If on the other hand $s^{2}\neq 1$, we have $s^{2p}=1$ and
    $\sum_{k=0}^{p-1}s^{2k}=0$ which yields
  \begin{equation*}
    J^{\Phi_{\sigma}}=
    \begin{pmatrix}
      s^{-2}t^{-1}(t-1)^{2}(s^{2}-1) & s^{-2}t^{-1}(1-t)(s^{2}-1)^{2}\\
      qt^{-1}(t-1)^{2} & q t^{-1}(1-t)(s^{2}-1)
    \end{pmatrix}.
  \end{equation*}
  As a result, $\Delta^{\Phi_{\sigma}}(t)=t^{-1}(t-1)\doteq t-1$ in this case.
\end{proof}

We can now proceed to prove the main result:%
\begin{proof}[Proof of Proposition~\ref{thm:Wpq-nonab-red-lie-lie-in-intersection}]
  Let $\prho\in\PRV(\Wpq)$ be non-abelian reducible.
  Lem\-ma~\ref{lem:nonab-red-specific-form-Wpq} shows that $s^{2}\neq 1$ implies
  $u^{2}=1$ (i.e.\ $u^{2}$ is the root of the polynomial $t-1$), and $s^{2}=1$
  implies $u^{2}$ is a root of $qt^{2}+(p-2q)t+q$. But these are precisely the
  twisted Alexander polynomials found in Lemma~\ref{lem:tw-alex-nonab-red-Wpq}
  for these cases, and so we have that $u^{2}$ is a simple root of the twisted
  Alexander polynomial, unless $s^{2}=1$ and $p/q=4$ (in which case $u^{2}=-1$
  is a double root).  Apart from this exception,
  Theorems~1.2~and~1.3 of \cite{MR2171800} (with $\delta=\mu_{1}$) then show
  that the character of the diagonal representation $\prho_{\alpha}$ (which is
  equal to the character of $\prho$) satisfies the claim of the Proposition.
\end{proof}


\subsection{A partially diagonal parametrization of p-reps}
\label{sec:param-irred-repr}

We determine a parametrization of p-reps similar to the normal form
\eqref{eq:normal-form-generators-RVC} from
section~\ref{sec:detour-through-eigenv}, but which is better suited to the
following cohomology calculations.

Let $\rho\in \RV(\Wpq)\subset\RV(W)$ be a p-rep. We therefore can assume that
$\rho(\mu_{0})$ is not parabolic, and that furthermore
$\trace \rho(\mu_{1})=\pm 2$.

\emph{For clarity of exposition, we will in the following restrict ourselves
  to the case $\trace \rho(\mu_{1})=+2$, and we attach the symbol ``$+$'' to
  all related sets and expressions to indicate this choice of sign. The
  negative trace case proceeds completely analogously and will be omitted
  here.}

As a consequence, $\rho$ can be conjugated into the form
\begin{equation}
  \label{eq:im-rho-diagonal}
  \rho(\mu_{0})=
  \begin{pmatrix}
    s & 0\\
    0 & s^{-1}
  \end{pmatrix}\qquad\text{and}\qquad
  \rho(\mu_{1})=
  \begin{pmatrix}
    a & -(a- 1)^{2}\\
    1 & 2-a
  \end{pmatrix},
\end{equation}
where $s\in \C^{\ast}$, $a\in \C$, subject to the condition (obtained from
applying \eqref{eq:im-rho-diagonal} to the first group relation)
\begin{equation}
  \label{eq:def-poly-r1}
  r_{1}^{+} = (a-1)\left[-(s^{2}-1)^2 a^2 + (s^{2}-1)(s^2-3) a-2\right]=0.
\end{equation}
Note that this implies $a\neq 0$ and $a\neq 2$ (otherwise $r_{1}^{+}=2s^{4}\neq 0$).

We will in the following use this condition to simplify all arising
polynomials modulo $r_{1}^{+}$, i.e.\ we will replace each polynomial $f$ by
its canonical form $f_{0}$ modulo the ideal generated by $r_{1}^{+}$, in
notation $f\equiv_{r_{1}^{+}} f_{0}$ (see e.g.\ \cite{MR1417938}). We also
extend this to rational functions by setting $f/g\equiv_{r_{1}^{+}} f_{0}/g$
if $f \equiv_{r_{1}^{+}} f_{0}$.

Using this parametrization, we then obtain for the longitude that
$\rho(\lambda_{0})=\begin{smallpmatrix} t & 0\\ 0 & t^{-1}\end{smallpmatrix}$
with $t\equiv_{r_{1}}a/(2-a)$, which implies that $\rho$ is diagonal on the
whole peripheral subgroup of the first boundary torus. This fact will be of great
importance in section~\ref{sec:cohomology}.

Applying \eqref{eq:im-rho-diagonal} to the filling relation,
$\rho(\mu_{0}^{p}\lambda_{0}^{q})=I$, now yields
\begin{align}
  \label{eq:def-poly-r2}
  r_{2}^{+}\equiv_{r_{1}^{+}}s^{p}\left(\tfrac{a}{2-a}\right)^{\!q}-1=0.
\end{align}
Note that $t\neq0$ implies $a\neq 0$ and $a\neq 2$.

Denote by $\RVD^{+}(\Wpq)\subset\RV(\Wpq)$ the subset of representations in
partially-diagonal normal form \eqref{eq:im-rho-diagonal} satisfying
$r_{1}^{+}=r_{2}^{+}=0$. 

Now note that $\rho\in\RVD^{+}(\Wpq)$ is reducible iff $a=1$. In this case the
polynomial $r_{1}^{+}$ is identically zero, and poses therefore no further
condition on $s$ and $a$. Then $\rho\in\RVC(\Wpq)\intersect\RVD^{+}(\Wpq)$,
and $a=u=1$ is equal to the eigenvalue of $\rho(\mu_{1})$. The condition
$r_{2}^{+}=0$ then implies $s^{p}=1$, just as we found in
section~\ref{sec:finding-preps}.

If $\rho\in\RVD^{+}(\Wpq)$ is irreducible, we have $a\neq 1$, and we can
further simplify $r_{1}^{+}$ to:%
\begin{equation}
  r_{1}^{+} = -(s^{2}-1)^2 a^2 + (s^{2}-1)(s^2-3) a-2.
\end{equation}
For simplicity, we will denote this polynomial again by $r_{1}^{+}$, since
there will be no danger of confusion.  Now note that the system
$r_{1}^{+}=r_{2}^{+}=0$ (for the coefficients $s$ and $a$) is completely
equivalent to the system $k_{1}=k_{2}=0$ (for eigenvalues $s$ and $t$) for
p-reps (with positive trace on $\rho(\mu_{1})$, i.e.\ $u=+1$) in
section~\ref{sec:finding-preps}, where we showed that its solutions are
characterized by the roots of the polynomial $\res_{p,q}$. Instead of
explicitly solving this system in the new parametrization, we will therefore
use the properties of $\res_{p,q}$ where appropriate (in particular in section
\ref{sec:deform-repr-degr}).


\subsection{The twisted cohomology and smoothness at p-reps}
\label{sec:cohomology}

This section is dedicated to proving
Proposition~\ref{thm:char-of-prep-is-smooth-point} for irreducible p-reps by
determining an explicit presentation for the cohomology group
$H^{1}:=H^{1}(\Wpq; \sltwist)$ $=H^{1}(\pi_{1}(\Wpq);\sltwist)$. Good background
references are \cite{MR672956} and \cite{MR1396960}, we also recommend the
articles \cite{MR1807271} and \cite{MR2171800} for a good summary in our
context.

\subsubsection{Preliminaries}

Consider the Lie group $\SL_{2}(\C)$ and its Lie algebra $\sll_{2}(\C)$, and
denote the adjoint homomorphism by $\Ad\map \SL_{2}(\C)\ra \Aut(\sll_{2}(\C))$
acting by conjugation, i.e.\ for $A\in\SL_{2}(\C)$ and $U\in\sll_{2}(\C)$ we
have $\Ad(A)(U)=AUA^{-1}$.

Let $\Gamma$ be a finitely presented group, and let
$\rho\in\RV(\Gamma)=\Hom(\Gamma;\SL_{2}(\C))$. The Lie algebra $\sll_{2}(\C)$
can then be viewed as a (left) $\Gamma$ module via the map
$\Ad\after\rho$, i.e.\ for $\gamma\in\Gamma$ and $A\in \SL_{2}(\C)$
we have the left action $\gamma\cdot A:= \Ad(\rho(\gamma))(A)$. Denote
$\sll_{2}(\C)$ with this left $\Gamma$ module structure by $\sltwist$.

Denote by $C^{\ast}(\Gamma;\sltwist)$ the space of
cochains with coefficients in $\sltwist$. Similarly, let
$B^{\ast}(\Gamma;\sltwist)$ and $Z^{\ast}(\Gamma;\sltwist)$ denote the
coboundaries and cocycles, respectively, and $H^{\ast}(\Gamma;\sltwist)$ the
cohomology group with coefficients in $\sltwist$.

Note that a cocycle $d\in Z^{1}(\Gamma;\sltwist)$ is a map $d\map \Gamma\ra
\sltwist$ satisfying
\begin{equation}
  \label{eq:cohomol-cocycle-condition}
  d(\gamma_{1}\gamma_{2})=d(\gamma_{1})+ \gamma_{1}\cdot d(\gamma_{2})
  \qquad \text{for all }\gamma_{1},\gamma_{2}\in\Gamma,
\end{equation}
using the short notation for the $\Ad\after\rho$ action of $\Gamma$ on $\sll_{2}(\C)$.
A cochain $b\map \Gamma\ra \sltwist$ is a coboundary if there exists an
$U\in\sll_{2}(\C)$ such that $b(\gamma)=\gamma\cdot U-U$ for all $\gamma\in\Gamma$.

It was an observation of Weil (\cite{MR0169956}, cf.\ also \cite{MR1396960})
that the Zariski tangent space\index{Zariski!tangent space}
$T_{\rho}^{\Zar}(\RV(\Gamma))$ at a point $\rho\in \RV(\Gamma)$ can be
identified with a subspace of the space of $1$-cocycles $Z^{1}(G;\sltwist)$,
and thus $\dim T_{\rho}^{\Zar}(\RV(\Gamma))\leq \dim Z^{1}(\Gamma;\sltwist)$.
The knowledge of the dimension of the space of cocycles can therefore be used
to determine whether a given representation or its character is a smooth point
in the representation or character variety, respectively.

\subsubsection{Coboundaries at p-reps}

First we determine the coboundaries $B^{1}:=B^{1}(\Wpq; \sltwist)$ for a p-rep
in partially-diagonal normal form $\rho\in\RVD^{\pm}(\Wpq)$. As in
section~\ref{sec:param-irred-repr}, we will restrict ourselves
to the case $\rho\in\RVD^{+}$, as the negative trace case proceeds
analogously.

Let 
\begin{equation}
e_{1}=\begin{pmatrix}0&1\\0&0\end{pmatrix},\quad 
e_{2}=\begin{pmatrix}1&0\\0&-1\end{pmatrix}\quad\text{and}\quad
e_{3}=\begin{pmatrix}0&0\\1&0\end{pmatrix}
\end{equation}
be a basis for the Lie algebra $\sll_{2}:=\sll_{2}(\C)$. For a matrix $A\in
\sll_{2}$, let $u_{A}\in B^{1}$ be the $1$-coboundary
defined by $u_{A}(\gamma)=A-\gamma \cdot A$ for $\gamma\in
\pi_{1}(\Wpq)$. Given a representation $\rho\in\RVD^{+}(\Wpq)$, the coboundaries
are then generated by
\begin{align}
  u_{e_{1}}&: & \mu_{0}&\mapsto 
            \begin{pmatrix}0&1-s^{2}\\0&0\end{pmatrix}, & \mu_{1}&\mapsto
            \begin{pmatrix}a&1-a^{2}\\1&-a\end{pmatrix},\\
  u_{e_{2}}&: & \mu_{0}&\mapsto 
            \begin{pmatrix}0&0\\0&0\end{pmatrix}, & \mu_{1}&\mapsto
            \begin{pmatrix}2(a-1)^{2}&-2a(a-1)^{2}\\ 2(a-2)&-2(a-1)^{2}\end{pmatrix},\\
  u_{e_{3}}&: & \mu_{0}&\mapsto 
            \begin{pmatrix}0&0\\1-s^{-2}&0\end{pmatrix}, & \mu_{1}&\mapsto
            \begin{pmatrix}(2-a)(a-1)^{2}&(a-1)^{4}\\-(a-1)(a-3)&-(2-a)(a-1)^{2}\end{pmatrix}.
\end{align}

Since cochains are determined by their images on the generators, we will
implicitly use the isomorphism $\tau: \sll_{2}^{2}\ra \C^{6}$ given by
\begin{equation}
  \label{eq:isom-C6}
  \left( 
    \begin{pmatrix}
      x_{1} & x_{2}\\
      x_{3} & -x_{1}
    \end{pmatrix},
    \begin{pmatrix}
      y_{1} & y_{2}\\
      y_{3} & -y_{1}
    \end{pmatrix} \right) \mapsto (x_{1},x_{2},x_{3},y_{1},y_{2},y_{3})
\end{equation}
to identify a cochain $u$ with an element of $\C^{6}$. We can therefore think
of $Z^{1}:=Z^{1}(\Wpq;\sltwist)$ and $B^{1}$ as subspaces of $\C^{6}$.

All coboundaries are now of the form $\alpha_{1} u_{e_{1}}+\alpha_{2}
u_{e_{2}}+\alpha_{3} u_{e_{3}}$ with $\alpha_{i}\in \C$. Seen as elements of
$\C^{6}$, they hence lie in the span of the row vectors of the matrix
\begin{equation}
  \label{eq:coboundaries-in-matrix-form}
  \cB:=\begin{pmatrix}
    0& 1-s^2& 0& a& 1-a^2& 1\\
    0& 0& 0& 2(a-1)^{2}& -2a(a-1)^{2} & 2(a-2)\\
    0& 0& 1-s^{-2}& (2-a)(a-1)^{2}& (a-1)^{4}& -(a-1)(a-3)
  \end{pmatrix}.
\end{equation}
Since $s\neq\pm1$, one sees that the rank of this
matrix, and hence the dimension of the subspace of coboundaries, is $3$.

Note
that this remains true when $\rho$ is reducible, as then $a=1$ and $\cB$ becomes
\begin{equation}
  \label{eq:coboundaries-in-matrix-form-reducible}
  \cB^{\text{red}}=\begin{pmatrix}
    0& 1-s^2& 0& 1& 0& 1\\
    0& 0& 0& 0& 0 & -2\\
    0& 0& 1-s^{-2}& 0& 0& 0
  \end{pmatrix}.
\end{equation}

\subsubsection{Cocycle condition from the first relation}

Let $(B^{1})^{\perp}$ be the
orthogonal complement of the subspace of coboundaries $B^{1}$ in $\C^{6}$ with
respect to the standard inner product.  For each cohomology element
$\bar{u}\in H^{1}$ there then is a unique representative $u\in
(B^{1})^{\perp}$ such that $\bar{u}=u+B^{1}$. The orthogonal complement is
given by the kernel of the matrix $\cB$, and we can therefore identify $H^{1}$
with a subspace of $(B^{1})^{\perp}$ in $\C^{6}$. In particular, we can obtain
it as the kernel of the linear map given by the matrix $\cB$ adjoined with
additional rows corresponding to the cocycle versions of the group relations.

We therefore examine the effect of the first group relation on cocycles: Using
$w_{1}:=\mu_{1}\mu_{0}\mu_{1}\mu_{0}^{-1}\mu_{1}^{-1}\mu_{0}^{-1}\mu_{1}\mu_{0}
=
\mu_{0}\mu_{1}\mu_{0}^{-1}\mu_{1}^{-1}\mu_{0}^{-1}\mu_{1}\mu_{0}\mu_{1}=:w_{2}$,
we see that each cocycle $u\in Z^{1}$ has to satisfy the matrix equation
$u(w_{1})-u(w_{2})=0$. Assuming that $u(\mu_{0})=
\begin{smallpmatrix}
  x_{1} & x_{2}\\
  x_{3} & -x_{1}
\end{smallpmatrix}$ and $u(\mu_{1})=
\begin{smallpmatrix}
  y_{1} & y_{2}\\
  y_{3} & -y_{1}
\end{smallpmatrix}$ as in \eqref{eq:isom-C6}, we can iteratively use the
cocycle condition~\eqref{eq:cohomol-cocycle-condition} to expand this
expression, yielding a matrix with entries which are polynomials in $x_{1}$,
$x_{2}$, $x_{3}$, $y_{1}$, $y_{2}$ and $y_{3}$ with coefficients in $\C[a^{\pm
  1},s^{\pm 1}]$.

When $\rho$ is reducible (but non-abelian), we find that only the lower left
entry of $u(w_{1})=u(w_{2})$ gives a non-vanishing equation, which
we write in vector form as
\begin{equation}
  \label{eq:cocycle-relation-r1-reducible}
{\footnotesize  \begin{pmatrix}
    0\\ -2(s^{2}-1) \\ 0 \\ -2(s^{2}-1)^{2} \\ s^{-2}(s^{2}-1)^{2}(s^{4}-s^{2}-1) \\ 0
  \end{pmatrix}^{T}
  \begin{pmatrix}
    x_{1}\\x_{2}\\x_{3}\\y_{1}\\y_{2}\\y_{3}
  \end{pmatrix}= 0.}
\end{equation}
This vector can now easily be seen to be linearly independent from the row
vectors in the matrix $\cB^{\text{red}}$ in~\eqref{eq:coboundaries-in-matrix-form-reducible}.

On the other hand, when $\rho$ is irreducible (and hence $a\neq 1$), the
matrix obtained is far more complicated. After
dividing the entries by a common nonzero factor of
$s^{-8}(s^{2}-1)(a-1)$ and reducing all coefficient polynomials in $\C[a^{\pm
  1},s^{\pm 1}]$ modulo $r_{1}^{+}$, one finds that the upper left entry
yields the equation
\begin{equation}
  \label{eq:cocycle-vector-form-equation-template-1}
  v^{T}  \begin{pmatrix}
    x_{1}&x_{2}&x_{3}&y_{1}&y_{2}&y_{3}
  \end{pmatrix}^{T}= 0,
\end{equation}
where the vector $v$ is given by
\begin{equation}
  \label{eq:cocycle-relation-r1-irred}
{\footnotesize  \begin{pmatrix}
-4(a^3s^2-a^3-7a^2s^2+7a^2+2as^4+8as^2-14a+4s^2+8)\\
-2(a^2s^2-a^2-2as^4-as^2+3a+s^4-2s^2-2)\\
2(a^4s^2-a^4-9a^3s^2+9a^3+24a^2s^2-29a^2-as^6-16as^2+39a+s^6-2s^4-8s^2-18)\\
2s^{2}(-6as^2+3a^2s^2-as^4-a^2+3a-2)\\
as^6+2as^4-2as^2+a-2s^4+4s^2-2\\
-3a^3s^2+a^3+16a^2s^2-8a^2-4as^4-19as^2+17a-2s^2-10
  \end{pmatrix},}
\end{equation}
and $v^{T}$ denotes the transposed vector.  Denote the first coefficient of
$v$ (corresponding to $x_{1}$) by $\psi_{1}$. One way to see that $v$ is
linearly independent from the row vectors in the matrix $\cB$ in
\eqref{eq:coboundaries-in-matrix-form} is to verify that $\psi_{1}$ does not
vanish at a representation $\rho\in\RVD^{+}(\Wpq)$: Using the computer algebra
system Maple\texttrademark\ (\cite{maple}) to calculate a reduced Groebner basis with respect to the pure
lexicographic order $a>s$ for the ideal generated by $\psi_{1}$ and
$r_{1}^{+}$, we find that its first generator is a power of $s$, implying that
the only common zeros with $r_{1}^{+}$ lie in the subspace defined by $s=0$, a
contradiction to our assumptions.

We remark that the other entries of
$u(w_{1})-u(w_{2})$ yield equations which can be shown to be consequences of
the coboundary conditions and equation \eqref{eq:cocycle-relation-r1-irred}.

\subsubsection{Cocycle condition from the filling relation}

We now study the condition forced on the cocycles after applying
$\pqslope$-filling on the first cusp. Letting $u(\mu_{0})=
\begin{smallpmatrix}
  x_{1} & x_{2}\\
  x_{3} & -x_{1}
\end{smallpmatrix}$ and $u(\mu_{1})=
\begin{smallpmatrix}
  y_{1} & y_{2}\\
  y_{3} & -y_{1}
\end{smallpmatrix}$ as before, we first determine the image
  $u(\lambda_{0})=u(\mu_{1}\mu_{0}\mu_{1}^{-1}\mu_{0}^{-1}\mu_{1}^{-1}\mu_{0}\mu_{1}\mu_{0}^{-1})$.
  A direct calculation shows that
\begin{equation}
  \label{eq:cocycle-image-longitude0}
  u(\lambda_{0})=:
\begin{pmatrix}
  l_{1} & l_{2}\\
  l_{3} & -l_{1}
\end{pmatrix},
\end{equation}
where the normalized entries $l_{1}$, $l_{2}$ and $l_{3}$ are given by
\begin{equation}
  \label{ex:cocycle-image-longitude0-details-l1}
{\footnotesize
  \begin{aligned}
    l_{1}=&-4s^{-6}(a^3s^2-a^3-2a^2s^2+4a^2-as^6-2as^4-as^2-5a+s^6+2s^4+2s^2+2)\\
    &-2s^{-8}(a^2s^2-a^2-as^2+3a-2s^6-2s^4-2s^2-2)\\
    &\begin{gathered}[t]
      +2s^{-4}(a^4s^2-a^4-3a^3s^2+5a^3-3a^2s^2-7a^2+2as^4\\+11as^2+a-2s^4-2s^2+2)
    \end{gathered}\\
    &-4s^{-2}(a^3s^2-a^3-3a^2s^2+4a^2+2as^2-5a+2)\\
    &-s^{-4}(2a^2s^2-2a^2-as^6-3as^2+6a-4s^4-2s^2-4)\\
    &\begin{gathered}[t]
      +2s^{-6}(a^4s^2-a^4-6a^3s^2+7a^3+12a^2s^2-19a^2\\-5as^2+23a-2s^4-6s^2-10),
    \end{gathered}
  \end{aligned}}
\end{equation}
\begin{equation}
  \label{ex:cocycle-image-longitude0-details-l2}
{\footnotesize
  \begin{aligned}
    l_{2}=
    &4s^{-6}(a^4s^2-a^4-a^3s^2+3a^3-11a^2s^2+3a^2+17as^2-15a+2s^2+10)\\
    &+s^{-8}(2a^3s^2-2a^3+4a^2-15as^4+a-s^4-8s^2-3)\\
    &\begin{gathered}[t]
      -s^{-4}(2a^5s^2-2a^5-4a^4s^2+8a^4-12a^3s^2-6a^3+32a^2s^2-12a^2\\
      -as^4-20as^2+21a+s^4-9)
    \end{gathered}\\
    &+2s^{-2}(2a^3s^2-2a^3-5a^2s^2+7a^2+2as^2-8a+s^2+3)a\\
    &+s^{-4}(2a^3s^2-2a^3-5a^2s^2+7a^2-5as^4+2as^2-7a+s^4-s^2+2)\\
    &\begin{gathered}[t]
      -s^{-6}(2a^5s^2-2a^5-11a^4s^2+13a^4+18a^3s^2-32a^3+4a^2s^2\\
      +30a^2-28as^2+2a-s^2-11)
    \end{gathered}
  \end{aligned}}
\end{equation}
and
\begin{equation}
  \label{ex:cocycle-image-longitude0-details-l3}
{\footnotesize
  \begin{aligned}
    l_{3}=
    &-4s^{-6}(a^2s^2-a^2-2as^8+2as^6-as^2+3a+2s^8-6s^6-2s^4-2s^2-2)\\
    &+s^{-2}(as^6-6as^4+7as^2-2a-s^6+8s^4-15s^2+4)\\
    &\begin{gathered}[t]
      +s^{-4}(2a^3s^2-2a^3-8a^2s^2+12a^2+as^{12}-6as^{10}+11as^8-4as^6-4as^4\\
      +2as^2-22a-s^{12}+8s^{10}-19s^8+12s^6+12s^4+12s^2+12)
    \end{gathered}\\
    &-2s^{-2}(2a^2s^2-2a^2-2as^4-3as^2+7a+2s^4-2s^2-6)\\
    &\begin{gathered}[t]
      +s^{-4}(as^{10}-7as^8+13as^6-12as^4+6as^2-a-s^{10}\\
      +9s^8-22s^6+22s^4-12s^2+2)
    \end{gathered}\\
    &\begin{gathered}[t]
      +s^{-6}(2a^3s^2-2a^3-13a^2s^2+15a^2+as^8+5as^6+11as^2\\
      -31a-s^8-3s^6+10s^4+14s^2+18).
    \end{gathered}
  \end{aligned}}
\end{equation}
The relation $\mu_{0}^{p}\lambda_{0}^{q}=1$ then translates to
$u(\mu_{0}^{p}\lambda_{0}^{q})=\smallzeromatrix$, where
\begin{equation}
  \label{eq:cocycle-relation-pq-filling}
  \begin{aligned}
 &u(\mu_{0}^{p}\lambda_{0}^{q})=u(\mu_{0}^{p})+\mu_{0}^{p}\cdot u(\lambda_{0}^{q}) = \\
&= \left(1+\mu_{0}+\cdots+\mu_{0}^{p-1}\right)\cdot u(\mu_{0})+ \mu_{0}^{p}
   \left(1+\lambda_{0}+\cdots+ \lambda_{0}^{q-1}\right)\cdot u(\lambda_{0}) \\
&=\begin{pmatrix}
  px_{1}+ql_{1} & x_{2}\sum_{k=0}^{p-1}s^{2k}+s^{2p}l_{2}\sum_{j=0}^{q-1}t^{2j}\\
  x_{3}\sum_{k=0}^{p-1}s^{-2k}+s^{-2p}l_{3}\sum_{j=0}^{q-1}t^{-2j} & -px_{1}-ql_{1} 
  \end{pmatrix}.
\end{aligned}
\end{equation}
We remark that this computation is the reason for our choice of the
parametrization in section~\ref{sec:param-irred-repr}, as we could otherwise
not have arrived at a general form in terms of $p$ and $q$.

Now $s\neq \pm 1$ implies that we can use
\begin{align}
  \sum_{k=0}^{p-1}s^{2k}&=\frac{s^{2p}-1}{s^{2}-1}\\
\intertext{and}
\sum_{k=0}^{p-1}s^{-2k}&=s^{-2(p-1)}\sum_{k=0}^{p-1}s^{2k}=s^{-2(p-1)}\frac{s^{2p}-1}{s^{2}-1}.
\end{align}

We now need to distinguish between the reducible and
irreducible case: When $\rho$ is reducible, $a=1$ implies $t=1$ and hence
$\sum_{j=0}^{q-1}t^{2j}=q$. But then $s^{p}t^{q}=1$ also implies $s^{p}=1$, and
equation~\eqref{eq:cocycle-relation-pq-filling} becomes
\begin{equation}
  \label{eq:cocycle-relation-pq-filling-cont-reducible}
  \zeromatrix =\begin{pmatrix}
    px_{1}+ql_{1} & ql_{2}\\
    ql_{3} & -px_{1}-ql_{1}
  \end{pmatrix}.
\end{equation}
But note that $a=1$ implies that $l_{1}=(s^{2}-s^{-2})y_{2}$,
$l_{2}=0$ and $l_{3}=-2s^{-2}x_{2}+2(s^{-2}-1)y_{1}+(s^{-4}-2+s^{2})y_{2}$.
It is then easy to see that $ql_{3}=0$ is already a consequence of the
coboundary conditions, and so
the only new condition remaining
from~\eqref{eq:cocycle-relation-pq-filling-cont-reducible} is $(p/q)
x_{1}+l_{1}=0$. Adjoining this condition in row vector form together with the
vector corresponding to equation~\eqref{eq:cocycle-relation-r1-reducible} to
the matrix from
\eqref{eq:coboundaries-in-matrix-form-reducible}, we then obtain a
presentation matrix for the cohomology group $H^{1}(\Wpq; \sltwist)$, viewed
as a subspace of~$\C^{6}$:%
\begin{equation}\label{eq:pres-matrix-filled-Wpq-reducible}
{\footnotesize  \begin{pmatrix}
    0& 1-s^{2}& 0& 1& 0& 1\\
    0& 0      & 0& 0& 0& -2\\
    0& 0      &1-s^{2} &0 & 0 & 0\\
    0& 2(1-s^{2})& 0&-2(1-s^{2})^{2}& s^{-2}(s^{4}-s^{2}-1)(1-s^{2})^{2}& 0\\
    p/q & 0  & 0  & 0  & s^{-2}(s^{4}-1)  & 0  \\
  \end{pmatrix}.}
\end{equation}
The rank of this matrix can easily be seen to be $5$ for all $s\neq \pm 1$, and we
therefore obtain
\begin{equation}
  \label{eq:dim-cohomology-Wpq-reducible}
  \dim_{\C} H^{1}(\Wpq; \sltwist) = 1.
\end{equation}

On the other hand, when $\rho$ is irreducible, we can use $t\neq \pm 1$ to
rewrite \eqref{eq:cocycle-relation-pq-filling}~as
\begin{equation}
  \begin{split}
    u(\mu_{0}^{p}\lambda^{q})&=\begin{pmatrix}
  px_{1}+ql_{1} &
  \frac{s^{2p}-1}{s^{2}-1}x_{2}+s^{2p}\frac{t^{2q}-1}{t^{2}-1}l_{2}\\
  s^{-2(p-1)}\frac{s^{2p}-1}{s^{2}-1}x_{3}+s^{-2p}t^{-2(q-1)}\frac{t^{2q}-1}{t^{2}-1}l_{3}
  &  -px_{1}-ql_{1}
  \end{pmatrix}\\
&=\begin{pmatrix}
  px_{1}+ql_{1} &
  \frac{s^{2p}-1}{s^{2}-1}x_{2}+\frac{1-s^{2p}}{t^{2}-1}l_{2}\\
  s^{2}\frac{1-t^{2q}}{s^{2}-1}x_{3}-t^{2}\frac{1-t^{2q}}{t^{2}-1}l_{3}
  &  -px_{1}-ql_{1}
  \end{pmatrix}\\
&=\begin{pmatrix}
  px_{1}+ql_{1} &
  (s^{2p}-1)\left( \frac{1}{s^{2}-1}x_{2}-\frac{1}{t^{2}-1}l_{2} \right)\\
  s^{-2p}(s^{2p}-1) \left( \frac{s^{2}}{s^{2}-1}x_{3}-\frac{t^{2}}{t^{2}-1}l_{3} \right)
  &  -px_{1}-ql_{1}
  \end{pmatrix}.
\end{split}
\end{equation}
From the upper right and lower left entries of this matrix we now see that
either $s^{2p}=1$ (which we know cannot be true by
Lemma~\ref{lem:res-only-roots-on-unit-circle-are-pm1}, since $s$ is also a
root of $\res_{p,q}$), or $(t^{2}-1)x_{2}=(s^{2}-1)l_{2}$ and
$s^{2}(t^{2}-1)x_{3}=t^{2}(s^{2}-1)l_{3}$.  After some careful algebraic
manipulations, it can be verified that these last two conditions are in fact
consequences of the coboundary
conditions~\eqref{eq:coboundaries-in-matrix-form} and equation
\eqref{eq:cocycle-relation-r1-irred}. The only new information comes therefore
from the upper left matrix entry. Since $q>0$, we can rewrite this in the form
$l_{1}+(p/q)x_{1}=0$, which gives us (after normalizing with respect to
$r_{1}^{+}$) the condition
\begin{equation}
  \label{eq:cocycle-vector-form-equation-template-2}
  w^{T}  \begin{pmatrix}
    x_{1}&x_{2}&x_{3}&y_{1}&y_{2}&y_{3}
  \end{pmatrix}^{T}= 0,
\end{equation}
where the vector $w$ is given by
\begin{equation}
  \label{eq:px-qu-condition-irred}
  {\footnotesize  \begin{pmatrix}
      p/q-4s^{-6}(a^3s^2-a^3-2a^2s^2+4a^2-as^6-2as^4-as^2-5a+s^6+2s^4+2s^2+2)\\
      -2s^{-8}(a^2s^2-a^2-as^2+3a-2s^6-2s^4-2s^2-2)\\
      2s^{-4}(a^4s^2-a^4-3a^3s^2+5a^3-3a^2s^2-7a^2+2as^4+11as^2+a-2s^4-2s^2+2)\\
      -4s^{-2}(a^3s^2-a^3-3a^2s^2+4a^2+2as^2-5a+2)\\
      -s^{-4}(2a^2s^2-2a^2-as^6-3as^2+6a-4s^4-2s^2-4)\\
      2s^{-6}(a^4s^2-a^4-6a^3s^2+7a^3+12a^2s^2-19a^2-5as^2+23a-2s^4-6s^2-10)
    \end{pmatrix}.}
\end{equation}

Adjoining the two row vectors corresponding to equations
\eqref{eq:cocycle-relation-r1-irred} and \eqref{eq:px-qu-condition-irred} to
the matrix $B$, we now obtain a presentation matrix for
$H^{1}(\Wpq; \sltwist)$.  With some care and making sure no polynomial used in
division can be zero on $\RVD^{+}$, we can row-reduce this matrix to the form
\begin{equation}\label{eq:pres-matrix-filled-Wpq-sketch}
  \begin{pmatrix}
    \psi_{1} & 0  & 0  & 0  & \ast  & \ast  \\
    0& s^{2}-1& 0& 0& -1& -(a-1)^{-2}\\
    0& 0& s^{-2}-1& 0& s^{2}(a-1)^{2}& s^{2}\\
    0& 0& 0& (a-1)^{2}& -a(a-1)^{2}& a-2\\
    0& 0  & 0  & 0  & \ast   & \ast
  \end{pmatrix},
\end{equation}
where each $\ast$ denotes a polynomial in $a$ and $s$ (normalized modulo
$r_{1}^{+}$) which does not vanish on $\RVD^{+}$, and $\psi_{1}$ is the first
coefficient of the vector in equation~\eqref{eq:cocycle-relation-r1-irred}.
This shows that the rank of this matrix is $5$, and we therefore find that
\begin{align}
  \label{eq:dim-cohomology-Wpq-summary}
  \dim_{\C} H^{1}(\Wpq; \sltwist) &= 1
\intertext{or equivalently,}
  \dim_{\C} Z^{1}(\Wpq; \sltwist) &= 4.
\end{align}

Since Corollary~\ref{thm:Wpq-nonab-red-contribute-to-norm} showed
that a non-abelian reducible p-rep lies on a non-trivial curve in $\RV(\Wpq)$,
and an irreducible p-rep does so automatically, 
 Theorem A of \cite{MR1909000} now shows:%
\begin{Proposition}\label{thm:char-of-prep-smooth-point}
  The character $\chi_{\rho}$ of a (reducible or irreducible) p-prep
  $\rho\in\RV(\Wpq)$ is a smooth point of $\CVi(W)$ and is contained in a
  unique non-trivial curve. Furthermore, we have
  $T_{\chi_{\rho}}^{\Zar}(\CV(\Wpq))\isom H^{1}(\Wpq; \sltwist)$, presented as
  a subset of $\C^{6}$ as the kernel of the matrix
  \eqref{eq:pres-matrix-filled-Wpq-reducible} in the reducible, and the matrix
  \eqref{eq:pres-matrix-filled-Wpq-sketch} in the irreducible case. \qed
\end{Proposition}
This then concludes the proof of
Proposition~\ref{thm:char-of-prep-is-smooth-point}, and furthermore provides 
an explicit parametrization for non-trivial cocycles in
$Z^{1}(M;\sltwist)$ at a given representation $\rho\in\RVD^{+}$. This will be
needed in the following section, when we will study deformations of these
representations.

\subsection{Deformations of representations and degree of the trace function}
\label{sec:deform-repr-degr}

We now turn our attention to the proof of
Proposition~\ref{thm:char-of-prep-is-simple-zero}. 

Since Proposition~\ref{thm:char-of-prep-smooth-point} showed that $\CVi(\Wpq)$
is smooth and positive-dimensional at the character of a p-rep $\rho$, we can
use the existence of non-trivial cocycles to find a smooth deformation
$\rho_{\nu}$ of $\rho_{0}=\rho$, which admits an expansion of the form
\begin{equation}
  \label{eq:expansion-of-deformation-of-rep}
  \rho_{\nu}=\exp\left(\nu u_{1}+\nu^{2}u_{2}+\nu^{3}u_{3}+\dots\right)\rho
\end{equation}
where $u_{j}: \pi_{1}(G)\ra\sll_{2}(\C)$, and $u_{1}\in Z^{1}(G;
\sltwist)$ (see e.g.\ \cite[\S~3]{MR1807271} or \cite{MR762512}).

Developing $f_{\chi_{\rho_{\nu}}}=\left(\trace \rho_{\nu}(\gamma)
\right)^{2}-4$ using the expansion
\begin{equation}
  \exp(A)=\sum_{j=1}^{\infty}A^{j}/j!
\end{equation}
for $A\in\sll_{2}(\C)$, we find that the trace function
$f_{\mu_{1}}(\chi_{\rho_{\nu}})$ has the following Taylor expansion at
$\rho_{0}$ (see e.g.\ Lemma 1.3 in \cite{MR1829566}):%
\begin{equation}
  \label{ex:taylor-expansion-fmu1}
f_{\mu_{1}}(\chi_{\rho_{\nu}})=f_{\mu_{1}}(\chi_{\rho_{0}})+2I_{\mu_{1}}(\chi_{\rho_{0}})\trace\bigl(u_{1}(\mu_{1})\rho_{0}(\mu_{1})\bigr)\nu+O(\nu^{2}).
\end{equation}
Note that in our case $G:=\pi_{1}(\Wpq)$ and
$I_{\mu_{1}}(\chi_{\rho_{0}})=\left(\trace\rho_{0}(\mu_{1})\right)^{2}=4$ for
$\rho_{0}\in\fD^{\pm}$. Since $f_{\mu_{1}}(\chi_{\rho_{0}})=0$ at a p-rep
$\rho_{0}$, the order of this zero is $1$ iff there is a non-trivial cocycle $u_{1}\in
Z^{1}(\Wpq;\sltwist)$ satisfying
\begin{equation}
  \label{eq:first-term-Taylor-non-zero}
  \trace\bigl(u_{1}(\mu_{1})\rho_{0}(\mu_{1})\bigr)\neq 0.
\end{equation}

To prove this, note that for $\rho_{0}\in\RVD^{+}$ and an arbitrary cochain $u_{1}$ we have
\begin{equation}
  \label{eq:first-term-Taylor-details}
  \trace\bigl(u_{1}(\mu_{1})\rho_{0}(\mu_{1})\bigr)
  \begin{aligned}[t]
    &= \trace
    \begin{pmatrix}
      y_{1} & y_{2}\\
      y_{3} & -y_{1}
    \end{pmatrix}
    \begin{pmatrix}
      a & -(a - 1)^{2}\\
      1 & 2-a
    \end{pmatrix}\\
    &= 2(a - 1)y_{1}+ y_{2} -(a - 1)^{2}y_{3},
  \end{aligned}
\end{equation}
and an analogous expression for $\rho_{0}\in\RVD^{-}$.

Let $P$ be the matrix obtained by adjoining the row vector corresponding to
this expression to the presentation matrix for $H^{1}(\Wpq;\sltwist)$ from
\eqref{eq:pres-matrix-filled-Wpq-reducible} or
\eqref{eq:pres-matrix-filled-Wpq-sketch}, respectively.  The
existence of a non-trivial cocycle $u_{1}$ such that
\eqref{eq:first-term-Taylor-details} is equal to a non-zero value
$c\in\C^{\ast}$ is now equivalent to finding a solution
$u_{1}:=(x_{1},x_{2},x_{3},y_{1},y_{2},y_{3})$ to the linear equation system
\begin{equation}\label{eq:trace-non-zero-linear-system}
  P  \left( x_{1}, x_{2}, x_{3}, y_{1}, y_{2}, y_{3} \right)^{T}
  =  \left( 0 ,  0 ,  0 ,  0 ,  0 ,  c  \right)^{T}
\end{equation}
and hence a non-trivial element in $H^{1}(\Wpq)$ represented by a cocycle
$u_{1}$ which satisfies \eqref{eq:first-term-Taylor-non-zero}.
In particular, if we can show that $\det P\neq 0$, the
system~\eqref{eq:trace-non-zero-linear-system} will in fact have a \emph{unique}
solution $u_{1}$ for any given $c\in\C^{\ast}$. 

Again restricting attention to $\rho_{0}\in\RVD^{+}$, we will first discuss the
case when $\rho_{0}$ is reducible: In this case, the matrix $P$ is given by
\begin{equation}
  \label{eq:final-matrix-P-reducible}
  P={\footnotesize \begin{pmatrix}
    0& 1-s^{2}& 0& 1& 0& 1\\
    0& 0      & 0& 0& 0& -2\\
    0& 0      &1-s^{2} &0 & 0 & 0\\
    0& 2(1-s^{2})& 0&-2(1-s^{2})^{2}& s^{-2}(s^{4}-s^{2}-1)(1-s^{2})^{2}& 0\\
    p/q & 0  & 0  & 0  & s^{-2}(s^{4}-1)  & 0  \\
    0& 0      & 0& 0& 1& 0\\
  \end{pmatrix}},
\end{equation}
and we find $\det P=-(4 p/q) (1-s^{2})^{2} (s^{4}-2s^{2}+2)$.  Now
recall from section~\ref{sec:red-solut} that the eigenvalue $s$ of a reducible
p-rep must satisfy $s^{p}=1$, $s\neq \pm 1$. None of the roots of
$s^{4}-2s^{2}+2=0$ lie on the unit circle, and so $\det P\neq 0$ unless
$p=0$: In that case, the trace in \eqref{eq:first-term-Taylor-details} will indeed be
zero for all cocycles $u_{1}$, and hence the order of a zero of $f_{\mu_{1}}$
at the character of a reducible non-abelian p-rep is indeed greater than $1$.
This also reflects the fact that there is a whole curve of non-abelian
reducible p-reps for $s\in \C^{\ast}-\{\pm 1\}$ in this case (since $p=0$,
$q=1$ implies $s^{p}t^{q}=t=1$, and there are no further conditions on $s$).

This proves Proposition~\ref{thm:char-of-prep-is-simple-zero} in the case of
a non-abelian reducible p-rep, $p/q\neq 0$.

We will now discuss the remaining case, when $\rho_{0}\in\RVD^{+}$ is irreducible.
Appending the row vector corresponding to the expression
\eqref{eq:first-term-Taylor-details} to the matrix \eqref{eq:pres-matrix-filled-Wpq-sketch} and applying Gauss-Jordan to reduce
the matrix, we obtain the $6\times 6$-matrix
\begin{equation}
  \label{eq:final-matrix-sketch}
  P=\begin{pmatrix}
    \psi_{1} &      &      &      & \ast                   & \ast         \\
         & \ast &      &      & \ast                   & \ast         \\
         &      & \ast &      & \ast                   & \ast         \\
         &      &      & \ast & \ast                   & \ast         \\
         &      &      &      & \ast                   & \ast         \\
         &      &      &      & 2a^{2}-2a+1 & (1-a)^{-1}(a^{3}-3a^{2}+5a-5)  \\
  \end{pmatrix}
\end{equation}
Denote by $d$ the determinant of the lower right $2\times 2$-block of this
matrix. We wish to show that this determinant is not zero for
$\rho_{0}\in\RVD^{+}$.  To show this, we again used the computer algebra
system Maple\texttrademark\ to calculate a Groebner basis of the ideal
generated by $r_{1}^{+}$ and the numerator of $d$ with respect to a
lexicographic ordering with $a>s$ (we remark that this computation
took more than 24 hours on a modern PC). We find that the first basis element
is given by $s^{30} d_{1} d_{2}$ where
\begin{equation}
  \label{ex:first-basis-element-groebner-determinant-h1}
  d_{1}=p(p-4q)s^4+(-6 p^2+24pq-32 q^2)s^2+p(p-4q)
\end{equation}
and
\begin{multline}
  \label{ex:first-basis-element-groebner-determinant-h2}
  d_{2}= 22 s^{40}-438 s^{38}+4185 s^{36}-24868 s^{34}+101875 s^{32}-304088 s^{30}\\
  +683740 s^{28}-1182928 s^{26}+1598312 s^{24}-1708564 s^{22}+1466502 s^{20}\\
  -1027864s^{18}+595850s^{16}-286072s^{14}+115452s^{12}-41808s^{10}\\
  +14586s^8-4582s^6+1097s^4-164s^2+11.
\end{multline}
If $d$ were to vanish at an irreducible representation $\rho_{0}\in\RVD^{+}$, each of the
polynomials in the Groebner basis, and in particular the first element, would
have to vanish.  We will therefore show that neither $d_{1}$ nor $d_{2}$
vanish at a representation $\rho_{0}\in\RVD^{+}$, which shows that $d$ cannot
vanish there.  For this, we will make use of the fact that the parametrization
of any irreducible representation $\rho_{0}\in\RVD^{\pm}$ also satisfies $\res_{p,q}=0$
(as determined in section~\ref{sec:irred-solut}). It suffices
therefore to show that $\res_{p,q}$ has no roots in common with $d_{1}$ and
$d_{2}$.

As for $d_{1}$, we note that its four roots are
\begin{equation}
  \label{eq:roots-of-h1}
\pm\frac{\left(\left(3(p-2q)^{2}+4q^2\right)\pm\left(2|p-2q|\sqrt{2(p-2q)^{2}+8q^2}\right)\right)^{1/2}}{\sqrt{p(p-4q)}}
\end{equation}
It is easy to see that in fact these roots are all
real-valued when $p>4q>0$ or when $p<0<4q$, and all imaginary when $0<p<4q$.
To see this, simply note that
\begin{equation}
  \left(3(p-2q)^{2}+4q^2\right)^{2}-\left(2|p-2q|\sqrt{2(p-2q)^{2}+8q^2}\right)^{2}
  =p^{2}(p-4q)^{2}\geq 0
\end{equation}
which shows that the numerator in \eqref{eq:roots-of-h1} is always real-valued.

But note that for $p>4q>0$ or $p<0<4q$, apart from $\pm 1$, none of the
non-zero roots of $\res_{p,q}$ are real-valued by
Lemma~\ref{thm:res-real-roots}. Similarly for $0<p<4q$, none of its non-zero
roots are pure imaginary by Lemma~\ref{thm:res-imaginary-roots}. So $d_{1}$ in
fact does not share any roots with $\res_{p,q}$.

It now remains to be shown that the polynomial $d_{2}$ has no common roots
with $\res_{p,q}$. To see this, first observe that $d_{2}$ is an irreducible
polynomial, while its leading coefficient is $22\neq 1$. On the other hand,
while $\res_{p,q}$ is in general \emph{not} irreducible (consider for example
$q=1$ and $p$ even), it is a \emph{monic} polynomial with integer
coefficients --- and therefore so are its irreducible factors: The Chebyshev
polynomial $T_{q}(y)$ has leading term $2^{q-1}y^{q}$, and substituting
$y=-s^{2}/2+2-s^{-2}/2$ gives the leading term $1/2\, s^{2q}$. Now note that the
``outer'' terms of $\res_{p,q}$ are $s^{p-2q}$ and $s^{-p+2q}$, and the middle
terms are given by $2T_{q}(y)$.

We conclude that all roots of
$\res_{p,q}$ are \emph{algebraic integers}, while this is not true for any of
the roots of $d_{2}$. This shows that $d_{2}$ and $\res_{p,q}$ have no roots
in common. We show the $40$ roots of $d_{2}$ in
figure~\ref{fig:roots-to-be-avoided-h2}.
\begin{figure}[htbp]
  \centering
  \includegraphics[width=0.6\textwidth]{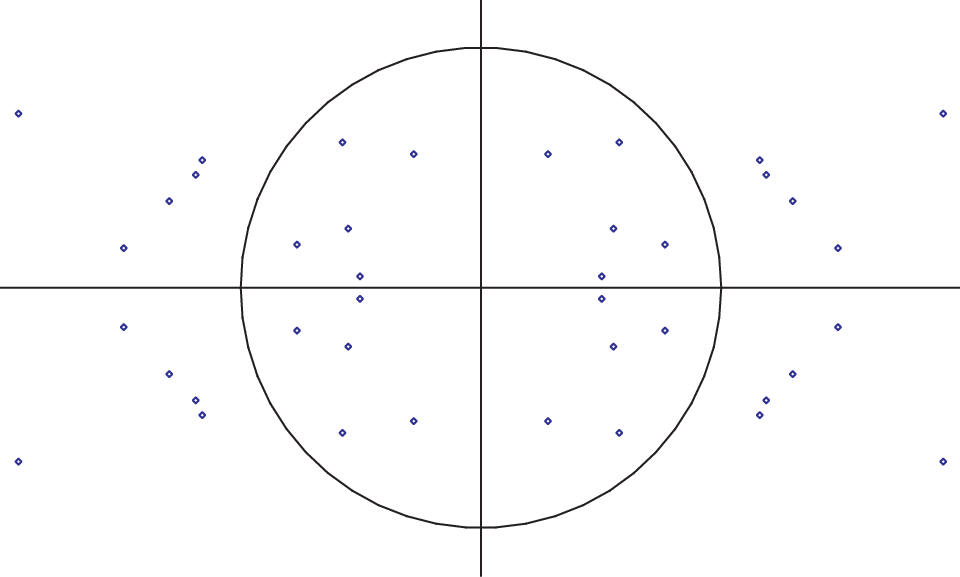}%
  \caption{The roots of the polynomial $d_{2}$}
  \label{fig:roots-to-be-avoided-h2}
\end{figure}

We can therefore conclude that the $2\times 2$ minor $d$ is indeed non-zero at
characters of p-reps, and that hence $\rank P=6$ and $\det
P\neq 0$, proving Proposition~\ref{thm:char-of-prep-is-simple-zero} for
irreducible p-reps.

This then concludes the proof of
Theorem~\ref{thm:minimal-seminorm-is-number-of-characters-of-preps}, showing
that the minimal seminorm $s$ is given by the number of characters of p-reps.



\section{Seifert fillings}
\label{sec:whitehead-seifert-fillings}

The goal of this section is to prove:

\begin{Proposition}\label{thm:norm-seifert-slopes}
  For $p$ odd, the values of the total Culler-Shalen seminorm of the 
  slopes $1$, $2$ and $3$ of $\Wpq$ are given by
  \begin{align*}
    \CSnorm{1}&=s + 2|p-6q|-2 \text{\quad (if\  $p/q\neq 6$)},\\
    \CSnorm{2}&=
      s+3|p-4q|-3  \text{\quad (if\  $p/q\neq 4$),} \\
    \CSnorm{3}&=s + 4|p-3q|-4 \text{\quad (if\  $p/q\neq 3$)},
  \end{align*}
  where $s$ is the minimal total Culler-Shalen seminorm.
\end{Proposition}
Note that while $p/q\neq 6$ and $p/q\neq 4$ are already satisfied when $p$ is
odd, we mention these conditions for completeness.

The proof is based on the well-known fact (see e.g.\ \cite[\S~9]{math.GT/0204228} or \cite{MR1988291}) that for
$\sigma\in\{1,2,3\}$, the closed manifold $\Wpq(\sigma)=W(\pqslope,\sigma)$ is
Seifert fibered with base orbifold a sphere with three cone
points (such manifolds are often called \emph{small Seifert}):%
\begin{itemize}
\item $W(\pqslope,1)$ is Seifert fibered over $S^{2}(2,3,|p-6q|)$ if $p/q\neq 6$
\item $W(\pqslope,2)$ is Seifert fibered over $S^{2}(2,4,|p-4q|)$ if $p/q\neq 4$
\item $W(\pqslope,3)$ is Seifert fibered over $S^{2}(3,3,|p-3q|)$ if $p/q\neq 3$
\end{itemize}
This can be seen from the fact that $W(1)$, $W(2)$ and $W(3)$ are all Seifert
fibred with base orbifold a disc with two cone points: $D^{2}(2,3)$,
$D^{2}(2,4)$ and $D^{2}(3,3)$, respectively. Furthermore, the corresponding
fibre slopes are $6$, $4$ and $3$, respectively. Therefore, applying a second
Dehn filling generically produces a Seifert fibred manifold with base orbifold
a sphere with three cone points, with only one exception each: When filling
along the fibre slope, the resulting closed manifold is a connected sum of lens spaces:
$W(1,6)\isom L_{2}\connsum L_{3}$, $W(2,4)\isom L_{2}\connsum L_{4}$ and
$W(3,3)\isom L_{3}\connsum L_{3}$.

The methods from \cite{MR1829566} now enable us to explictly determine
the Culler-Shalen seminorm of these small Seifert manifolds:

  Let $\sigma\in \{1,2,3\}$ be one of the Seifert filling slopes for $\Wpq$.
  Theorem~C of \cite{MR1829566} shows that in this situation
  $\CSnorm{\sigma}=s+2A$, where $A$ is the number of characters of non-abelian
  representations $\rho\in\RV(\Wpq)$ satisfying $\rho(\sigma)=\pm I$ whose
  characters lie on non-trivial curves in $\CV(\Wpq)$.
  
  Note that such a representation necessarily induces a (non-abelian)
  $\PSL_{2}(\C)$ representation $\prho\in\PRV(\Wpq)$ with
  $\prho(\sigma)=\pm I$, which in turn factors through to give a
  representation $\prho'\in\PRV(\Wpq(\sigma))$.  (An abelian irreducible
  representation $\prho\in\PRV$ necessarily has image isomorphic to
  $\Zmod{2}\oplus\Zmod{2}$, generated by
  $\pm \begin{smallpmatrix} i & 0 \\ 0 & -i \end{smallpmatrix}$ and
  $\pm \begin{smallpmatrix} 0 & -b \\ b^{-1} & 0 \end{smallpmatrix}$ for some
  $b\in\C^{\ast}$.  This possibility is excluded in our situation by the fact
  that $p$ is odd.)  Since $\Wpq(\sigma)$ has base orbifold $S^{2}(a,b,c)$, a
  sphere with three cone points, there is a presentation of
  $\pi_{1}(\Wpq(\sigma))$ of the form
  \begin{equation}
    \label{eq:presentation-small-Seifert-fund-group}
    \pi_{1}(\Wpq(\sigma))=\genrel[x,y,h]{h \text{ central},\ x^{a}=h^{a'},\ y^{b}=h^{b'},\
      (xy)^{c}=h^{c'}}
  \end{equation}
  where $\gcd(a,a')=\gcd(b,b')=\gcd(c,c')=1$ (cf.\ \cite[\S~VI]{MR565450}).

  Lemma~3.1 from \cite{MR1829566} now shows that the representation
  $\prho'\in\PRV(\Wpq(\sigma))$ again factors through to give a representation
  $\prho''\in\PRV(\Delta(a,b,c))$, where $\Delta(a,b,c)\isom
  \pi_{1}^{\text{orb}}(S^{2}(a,b,c))$ is the $(a,b,c)$ triangle group with
  presentation
  \begin{equation}
    \Delta(a,b,c)=\genrel[x,y]{x^{a}=y^{b}=(xy)^{c}=1}.
  \end{equation}
  Conversely, a given representation $\prho''\in\PRV(\Delta(a,b,c))$ induces a
  unique representation $\prho'\in\PRV(\Wpq(\sigma))$, which in turn induces a
  unique representation $\prho\in\PRV(\Wpq)$ satisfying $\prho(\sigma)=\pm I$.
  We summarize the situation in the following diagram:%
  \begin{equation}\label{eq:diagram-of-representations-Seifert}
    \xymatrix{
      \pi_{1}(\Wpq) \ar[r]^{\rho}\ar[d]\ar[dr]^{\prho} & \SL_{2}(\C)\ar[d]\\
      \pi_{1}(\Wpq(\sigma))\ar[d] \ar[r]^{\prho'} & \PSL_{2}(\C)\\
      \Delta(a,b,c)\ar[ur]^{\prho''} & 
    }
  \end{equation}

  Using the following theorem from \cite{MR1829566}, which we restate here
  in slightly adapted form, together with the fact that
  $H_{1}(\Wpq(r/s))=\Zmod{|p|}\oplus \Zmod{|r|}$, one can now explicitly
  determine the number of $\PSL_{2}(\C)$ characters of $\Delta(a,b,c)$, and
  hence of $\pi_{1}(\Wpq)$, for each of the Seifert
  filling slopes:
  
\begin{Theorem*}[\cite{MR1829566}, Theorem D]~\\
  Let $M$ be a non-Haken small Seifert manifold, with base orbifold of the
  form $S^{2}(a,b,c)$ where $a,b,c\geq 2$.
  \begin{enumerate}
  \item The number of reducible $\PSL_{2}(\C)$-characters of $\pi_{1}(M)$ is
    \[ \floor{\frac{\abs{H_{1}(M)}}{2}}+1+
    \begin{cases}
      0 & \text{if $\gcd(a,b,c)$ is odd,}\\
      1 & \text{if $\gcd(a,b,c)$ is even.}
    \end{cases}
. \]
  \item The number of $\PSL_{2}(\C)$-characters of $\pi_{1}(M)$
    representations with image a dihedral group of order at least $4$ is
    \[
    \sigma(b,c)\bfloor{\frac{a}{2}} +\sigma(a,c)\bfloor{\frac{b}{2}}
    +\sigma(a,b)\bfloor{\frac{c}{2}}-2\sigma(a,b)\sigma(a,c)\sigma(b,c)  \]
    where $\sigma(m,n)=\floor{\frac{\gcd(2,m,n)}{2}}$ is $1$ if $m$ and $n$
    are both even, and $0$ otherwise.
  \item The total number of $\PSL_{2}(\C)$-characters of $\pi_{1}(M)$ is 
    \begin{gather*}
      \bfloor{\frac{\abs{H_{1}(M)}}{2}}+1
      + \bfloor{\frac{a}{2}}\bfloor{\frac{b}{2}}\bfloor{\frac{c}{2}}
      + \bfloor{\frac{a-1}{2}}\bfloor{\frac{b-1}{2}}\bfloor{\frac{c-1}{2}}\\
      - \bfloor{\frac{\gcd(ab,ac,bc)}{2}} + \bfloor{\frac{\gcd(a,b)}{2}} 
      + \bfloor{\frac{\gcd(a,c)}{2}} + \bfloor{\frac{\gcd(b,c)}{2}}
    \end{gather*}
  \end{enumerate}
\qed
\end{Theorem*}

The results are summarized in tables~\ref{tab:num-PSL-characters-Wrh-pq-1},
\ref{tab:num-PSL-characters-Wrh-pq-2} and
\ref{tab:num-PSL-characters-Wrh-pq-3} on
page~\pageref{tab:num-PSL-characters-Wrh-pq-1}.  Also included are the numbers
of non-abelian reducible characters, which can be found by direct calculation.

  \begin{sidewaystable}[htbp]
    \centering
    \renewcommand{\arraystretch}{1.2}
    \begin{tabular}{C||C|CC|CC}
      \gcd(6,p) &\text{total} & \text{irreducible} & \text{dihedral} & \text{reducible} & \text{non-ab.\ red.}\\
      \hline\hline

      1 & \tfrac{1}{2}\left(|p|+|p-6q|\right) 
      & \tfrac{1}{2}\left(|p-6q|-1\right) & 0 
      & \tfrac{1}{2}\left(|p|+1\right)    & 0\\
      3 & \tfrac{1}{2}\left(|p|+|p-6q|\right) 
      & \tfrac{1}{2}\left(|p-6q|-1\right) & 0
      & \tfrac{1}{2}\left(|p|+1\right)    & 0\\
      \hline
      2 & \tfrac{1}{2}\left(|p|+|p-6q|\right)+1 
      & \tfrac{1}{2}|p-6q|   & 1
      & \tfrac{1}{2}|p|+1    & 0\\
      6 & \tfrac{1}{2}\left(|p|+|p-6q|\right)
      & \tfrac{1}{2}|p-6q|-1   & 1
      & \tfrac{1}{2}|p|+1      & 1\\
    \end{tabular}
    \\[2ex]\caption{Number of $\PSL_2(\C)$-characters of $\pi_1 \Wpq(1)$}
    \label{tab:num-PSL-characters-Wrh-pq-1}

    \vspace{1cm}

    \begin{tabular}{C||C|CC|CC}
      \gcd(4,p) &\text{total} &  \text{irreducible} & \text{dihedral} &  \text{reducible} & \text{non-ab.\ red.}\\
      \hline\hline
      1 & |p|+|p-4q|    
      & |p-4q|-1 & \tfrac{1}{2}\left(|p-4q|-1\right) 
      & |p|+1    & 0\\
      \hline
      2 & |p|+|p-4q|+2  
      & |p-4q|   & \tfrac{1}{2}|p-4q|+1 
      & |p|+2    & 0\\ 
      4 & |p|+|p-4q|+1  
      & |p-4q|-1 & \tfrac{1}{2}|p-4q|+1
      & |p|+2    & 1\\
    \end{tabular}
    \\[2ex]\caption{Number of $\PSL_2(\C)$-characters of $\pi_1 \Wpq(2)$}
    \label{tab:num-PSL-characters-Wrh-pq-2}

    \vspace{1cm}

    \begin{tabular}{C||C|CC|CC}
      \gcd(6,p) &\text{total} & \text{irreducible} & \text{dihedral} & \text{reducible} & \text{non-ab.\ red.}\\
      \hline\hline
      1 & \tfrac{3}{2}\left(|p|-1\right)+|p-3q|   
      & |p-3q|-1                          & 0 
      & \tfrac{3}{2}\left(|p|-1\right)+1  & 0\\
      3 & \tfrac{3}{2}\left(|p|-1\right)+|p-3q|-1 
      & |p-3q|-2                          & 0 
      & \tfrac{3}{2}\left(|p|-1\right)+1  & 1\\
      \hline
      2 & \tfrac{3}{2}|p|+|p-3q|   
      & |p-3q|-1                 & 0 
      & \tfrac{3}{2}|p|+1        & 0\\
      6 & \tfrac{3}{2}|p|+|p-3q|-1 
      & |p-3q|-2                 & 0 
      & \tfrac{3}{2}|p|+1        & 1\\
    \end{tabular}
    \\[2ex]\caption{Number of $\PSL_2(\C)$-characters of $\pi_1 \Wpq(3)$}
    \label{tab:num-PSL-characters-Wrh-pq-3}

    \vspace*{-0.5cm}
  \end{sidewaystable}

  \newpage
  It now remains to determine the number of non-abelian $\SL_{2}(\C)$
  characters, i.e.\ to determine which of the $\PSL_{2}(\C)$ representations
  $\prho\in\PRV(\Wpq)$ lift to representations $\rho\in\RV(\Wpq)$. There is in
  general an obstruction to this lifting, which is contained in the cohomology
  group
  \begin{equation}
    \label{eq:obstruction-to-lifting-H2}
    H^{2}(\Wpq;\Zmod{2})\isom 
    \begin{cases}
      0 & \text{if $p$ is odd,}\\
      \Zmod{2} & \text{if $p$ is even.}
    \end{cases}
  \end{equation}
  When $p$ is even, we therefore possibly do have
  an obstruction to lifting, and it is to be expected that some
  representations do lift, but not others, making it impossible to directly
  make a statement about the number of $\SL_{2}(\C)$ representations.

  On the other hand, there is \emph{no obstruction to lifting} when $p$ is
  odd, \emph{and we shall therefore restrict our attention to this case.}

  The number of lifts of a given $\PSL_{2}(\C)$ representation
  is generically given by the cardinality of
  \begin{equation}
    \label{eq:obstruction-to-lifting-H1}
    H^{1}(\Wpq;\Zmod{2})\isom 
    \begin{cases}
      \Zmod{2} & \text{if $p$ is odd,}\\
      \Zmod{2}\oplus \Zmod{2} & \text{if $p$ is even,}
    \end{cases}
  \end{equation}
  but dihedral representations only lift half as many times as given by
  this number (see e.g.\ Lemma 5.5 of~\cite{BZfinite} and the remarks
  following it for more details).

  When $p$ is odd, dihedral $\PSL_{2}(\C)$ characters are therefore covered by
  only one $\SL_{2}(\C)$ character, whereas the remaining non-abelian
  characters are covered by two, giving us the results listed in
  table~\ref{tab:whitehead-seifert-nonabelian-SL2C-p-odd}, where the total
  number of non-abelian $\SL_{2}(\C)$ characters shown in the last column
  gives the constant $A$ needed to conclude the proof of
  Proposition~\ref{thm:norm-seifert-slopes}.

  \begin{table}[htbp]
    \centering
{\small
    \begin{tabular}{C|CCC|C}
      & \text{irred.\ non-dihedral} & \text{dihedral} & \text{non-ab.\ red.} &
      \text{total}\\
      \hline
      \Wpq(1) & |p-6q|-1 & 0 & 0 & |p-6q|-1\\
      \Wpq(2) & |p-4q|-1 & \frac{1}{2}(|p-4q|-1)& 0 & \frac{3}{2}\left(|p-4q|-1\right)\\
      \Wpq(3) &    \begin{cases}
        2|p-3q|-2 & \text{if $3\ndivides p$}\\
        2|p-3q|-4 & \text{if $3\divides p$}\\
      \end{cases} & 0 & 
      \begin{cases}
        0 & \text{if $3\ndivides p$}\\
        2 & \text{if $3\divides p$}\\
      \end{cases} & 2\left(|p-3q|-1\right) \\
    \end{tabular}
    \\[2ex]\caption{Number of non-abelian $\SL_2(\C)$ characters for Seifert fillings
      of $\Wpq$ when $p$ is odd}
    \label{tab:whitehead-seifert-nonabelian-SL2C-p-odd}
}  \end{table}



\section{The total Culler-Shalen seminorm}
\label{sec:total-norm-solving-systems}

Putting together the results of the previous sections, we now completely
determine the total Culler-Shalen seminorm of the manifolds $\Wpq$ for $p$ odd
and $p/q\neq 3$, proving Theorem~\ref{thm:total-cs-norm-Wpq}. In particular,
we will show that the upper bounds established in
Theorem~\ref{thm:minimal-seminorm-is-number-of-characters-of-preps} (cf.\
section~\ref{sec:minimal-CS-norm}) for the \emph{minimal} total Culler-Shalen
seminorm $s$ are in fact attained, which in turn shows Corollary~\ref{cor:roots-of-res-are-simple}.

Using the possible boundary slopes $\beta_{j}$ from
section~\ref{sec:boundary-slopes-Wpq}, we 
apply equation~\eqref{eq:total-CSnorm-as-linear-comb-of-dist} with unknown
coefficients $a_{j}$ and $s$ to the Seifert filling slopes $1$, $2$, $3$ and
compare with the explicit values (relative to $s$) obtained in
Proposition~\ref{thm:norm-seifert-slopes}.

From these equations, we then obtain a linear system which we can try to solve
for the unknown coefficients $a_{j}$ and $s$.
Note that depending on the range of the filling slope $\pqslope$, the
parametrization of the boundary slope $\beta_{2}$ changes (cf.\ 
table~\ref{tab:boundary-slopes-Wrh} on page
\pageref{tab:boundary-slopes-Wrh}), and hence we get different
distances and accordingly different linear systems for each range. We will
therefore discuss each range in sequence.

While there are some cases (when $p/q\in(3,4)$ or $p/q\in(4,6)$) when the
linear system already determines a unique solution, the system is in general
underdetermined and we need more information. This is where we use the upper
bound for the total minimal seminorm $s$ established in
Theorem~\ref{thm:minimal-seminorm-is-number-of-characters-of-preps}. In two sub-cases ($\pqslope\in (2,3)$ and
$\pqslope>6$), it turns out that even this is not yet sufficient, but using a
symmetry argument about the resultant $\res_{p,q}$, we can reduce the case to
one which was already proven before.

As before, we will always assume $q>0$ and $p\in \Z$ with $(p,q)=1$ and (due
to the restrictions from section~\ref{sec:whitehead-seifert-fillings}) $p$
odd as well as $p/q\neq 3$.


\subsection{\texorpdfstring{Assuming $\pqslope\in (-\infty,0)$, $p$ odd}%
  {Assuming p/q < 0, p odd}}
\label{sec:CS-negative-odd-filling}

In this case, we obtain the following distances between the boundary, Seifert
and meridional slopes:
\begin{table}[h!]
  \centering
  \begin{tabular}{C|CCCCCCC}
    \Delta & \beta_{1}=4 & \beta_{2}=\frac{4q}{p} & \beta_{3}=0 & 1 & 2 & 3 & \mu=\infty\\\hline
    \beta_{1} & 0 & 4(-p+q) & 4    & 3 & 2 & 1 & 1\\
    \beta_{2} &   &   0    & 4q   & -p+4q & -2p+4q & -3p+4q & -p\\
    \beta_{3} &   &        & 0    & 1 & 2 & 3 & 1
  \end{tabular}
  \\[2ex]\caption{Distances between slopes in $\Wpq$ for $\pqslope\in
    (-\infty,0)$, $p$ odd, $q>0$}
  \label{tab:distances-negative-odd-filling}
\end{table}

Using the expressions of the seminorms for the Seifert slopes from
Proposition~\ref{thm:norm-seifert-slopes} in terms of the minimal seminorm $s$ on
one hand and the corresponding linear combinations with the slope distances on
the other, we get the following linear system to solve for $s$, $a_{1}$,
$a_{2}$ and $a_{3}$:%
\begin{subequations}\label{eq:CS-negative-odd}
  \begin{alignat}{2}
    s + 2(-p+6q-1) &= \CSnorm{1}      &&= 3a_{1} + (-p+4q)a_{2}  + a_{3}\\
    s + 3(-p+4q-1) &= \CSnorm{2}      &&= 2a_{1} + (-2p+4q)a_{2} + 2a_{3}\\
    s + 4(-p+3q-1) &= \CSnorm{3}      &&=  a_{1} + (-3p+4q)a_{2} + 3a_{3}\\
    s              &= \CSnorm{\infty} &&=  a_{1} - pa_{2}        + a_{3}
  \end{alignat}
\end{subequations}
which is equivalent in matrix form to
\begin{equation}
  \label{eq:CS-negative-odd-matrixform}
  \underbrace{\begin{pmatrix}
    3 & -p+4q & 1 & -1\\
    2 & -2p+4q & 2 & -1\\
    1 & -3p+4q & 3  & -1\\
    1 & -p & 1 & -1
  \end{pmatrix}}_{A}\begin{pmatrix}a_{1}\\a_{2}\\a_{3}\\s\end{pmatrix}=
  \begin{pmatrix}
    2(-p-1)+12q\\
    3(-p-1)+12q\\
    4(-p-1)+12q\\
    0
  \end{pmatrix}
\end{equation}
Since the second row of $A$ is equal to one half of the sum of the first and third
rows, the rank of $A$ is seen to be three, and the system is therefore
underdetermined. We therefore need more information to determine the solution.

Using the bound from table~\ref{tab:number-of-conjug-classes-preps-Wpq}, we have $s\leq -3p+4q-3$. In other words, we have
\begin{equation}
  s=-3p+4q-3-z,
\end{equation}
where $z=z(p,q)\geq 0$ is some \emph{positive}, integer-valued function in $p$ and $q$.

Using this expression for $s$, we can formally solve the
system~\eqref{eq:CS-negative-odd-matrixform} to obtain the following solution
in terms of $z$:\begin{equation}
  \label{eq:CS-positive-odd-guessed-solution}
  \begin{pmatrix}
    a_{1}\\a_{2}\\a_{3}\\s
  \end{pmatrix}=
  \begin{pmatrix}
    -p+2q-1-\frac{z}{2}\\
    2+\frac{z}{4q}\\
    2q-2+z\frac{p-2q}{4q}\\
    -3p+4q-3-z
  \end{pmatrix}.
\end{equation}
We can now use the fact that $a_{1}$, $a_{2}$, $a_{3}$ and $s$ all are
\emph{even}, \emph{non-negative} integers: The expression for $a_{2}$ then implies
that $z$ is divisible by $8q$, i.e.\ $z=8qz'$ for some non-negative
integer-valued $z'=z'(p,q)$. Using this in the expression for $a_{3}$ then
shows that $2q-2+2z'(p-2q)\geq 0$, which implies that
\begin{equation}
  z'(p-2q)\geq 1-q.
\end{equation}
Since $p<0$ and $q>0$ implies $p-2q<0$, we then get
\begin{equation}
  z'\leq \frac{1-q}{p-2q}=\frac{q-1}{2(q+|p|/2)}<1.
\end{equation}
But since $z'\in\Z_{+}$, we must conclude that in fact $z'=z=0$.

This proves Theorem~\ref{thm:total-cs-norm-Wpq} and
Corollary~\ref{cor:roots-of-res-are-simple} for $p/q<0$.


\subsection{\texorpdfstring{Assuming $\pqslope\in (0,2)$, $p$ odd}%
  {Assuming 0<p/q<2, p odd}}
\label{sec:CS-positive-0-2-odd-filling}

In this case, we obtain the following distances between the relevant slopes:
\begin{table}[h!]
  \centering
  \begin{tabular}{C|CCCCCCC}
    \Delta & \beta_{1}=4 & \beta_{2}=\frac{2p+4q}{p} & \beta_{3}=0 & 1 & 2 & 3 & \mu=\infty\\\hline
    \beta_{1} & 0 & -2p+4q & 4    & 3 & 2 & 1 & 1\\
    \beta_{2} &   &   0    & 2p+4q   & p+4q & 4q & -p+4q & p\\
    \beta_{3} &   &        & 0    & 1 & 2 & 3 & 1
  \end{tabular}
  \\[2ex]\caption{Distances between slopes in $\Wpq$ for $\pqslope\in
    (0,2)$, $p$ odd}
  \label{tab:distances-positive-0-2-odd-filling}
\end{table}

The corresponding linear system is in this case given by the equations:
\begin{subequations}\label{eq:CS-positive-0-2-odd}
  \begin{alignat}{2}
    s + 2(-p+6q-1) &= \CSnorm{1}      &&= 3a_{1} + (p+4q) a_{2} +   a_{3}\\
    s + 3(-p+4q-1) &= \CSnorm{2}      &&= 2a_{1} + 4q a_{2}     + 2 a_{3}\\
    s + 4(-p+3q-1) &= \CSnorm{3}      &&= a_{1}  + (-p+4q)a_{2} + 3 a_{3}\\
    s              &= \CSnorm{\infty} &&= a_{1}  + p a_{2}      + a_{3}
  \end{alignat}
\end{subequations}
Using the upper bound $s\leq p+4q-3$ from
table~\ref{tab:number-of-conjug-classes-preps-Wpq}, we can assume
\begin{equation}
  s=p+4q-3-z,
\end{equation}
where $z=z(p,q)\geq 0$ is some \emph{positive}, integer-valued function in $p$ and $q$.

Using this expression for $s$, we can formally solve the
system~\eqref{eq:CS-positive-0-2-odd} to obtain the following solution
in terms of $z$:\begin{equation}
  \label{eq:CS-positive-0-2-odd-guessed-solution}
  \begin{pmatrix}
    a_{1}\\a_{2}\\a_{3}\\s
  \end{pmatrix}=
  \begin{pmatrix}
    -p+2q-1-\frac{qz}{2q-p}\\
    2+\frac{z}{2(2q-p)}\\
    2q-2-\frac{z}{2}\\
    p+4q-3-z
  \end{pmatrix}.
\end{equation}
Now recall that $a_{1}$, $a_{2}$, $a_{3}$ and $s$ all are
\emph{even}, \emph{non-negative} integers. The expression for $a_{2}$ then
implies that $z$ is divisible by $4(2q-p)$, i.e.\ $z=4(2q-p)z'$ for some
non-negative integer-valued $z'=z'(p,q)$. Using this in the expression for
$a_{1}$ then shows that $-p+2q-1-4qz'\geq 0$, which implies that
\begin{equation}
  z'\leq \frac{-p+2q-1}{4q}=1-\frac{p+1+2q}{4q}<1,
\end{equation}
so in fact $z'=z=0$, proving Theorem~\ref{thm:total-cs-norm-Wpq} and
Corollary~\ref{cor:roots-of-res-are-simple} in the case $p/q\in(0,2)$.


\subsection{\texorpdfstring{Assuming $\pqslope\in (2,3)\cup (3,4)$, $p$ odd}%
  {Assuming 2<p/q<3 and 3<p/q<4, p odd}}
\label{sec:CS-positive-2-4-odd-filling}

In this case, we obtain the following distances between the relevant slopes:
\begin{table}[h!]
  \centering
  \begin{tabular}{C|CCCCCCC}
    \Delta & \beta_{1}=4 & \beta_{2}=\frac{-p+6q}{q} & \beta_{3}=0 & 1 & 2 & 3 & \mu=\infty\\\hline
    \beta_{1} & 0 & p-2q & 4    & 3 & 2 & 1 & 1\\
    \beta_{2} &   &   0    & -p+6q   & -p+5q & -p+4q & |p-3q| & q\\
    \beta_{3} &   &        & 0    & 1 & 2 & 3 & 1
  \end{tabular}
  \\[2ex]\caption{Distances between slopes in $\Wpq$ for $\pqslope\in
    (2,4)$, $p$ odd}
  \label{tab:distances-positive-2-4-odd-filling}
\end{table}

In this case, the corresponding linear system is:
\begin{subequations}\label{eq:CS-positive-2-4-odd}
  \begin{alignat}{2}
    s + 2(-p+6q-1)  &= \CSnorm{1}      &&= 3a_{1} + (-p+5q) a_{2} +  a_{3}\\
    s + 3(-p+4q-1)  &= \CSnorm{2}      &&= 2a_{1} + (-p+4q) a_{2} + 2a_{3}\\
    s + 4(|p-3q|-1) &= \CSnorm{3}      &&=  a_{1} + |p-3q|  a_{2} + 3a_{3}\\
    s               &= \CSnorm{\infty} &&=  a_{1} +       q a_{2} +  a_{3}
  \end{alignat}
\end{subequations}
If $\pqslope\in (3,4)$, this system already has rank $4$, and we can determine the
unique solution without further assumptions:\begin{equation}
  \label{eq:CS-positive-3-4-odd-exact-solution}
  \begin{pmatrix}
    a_{1}\\a_{2}\\a_{3}\\s
  \end{pmatrix}=
  \begin{pmatrix}
    p-2q-1\\
    4 \\
    2q-2\\
    p+4q-3
  \end{pmatrix}.
\end{equation}
This proves Theorem~\ref{thm:total-cs-norm-Wpq} and
Corollary~\ref{cor:roots-of-res-are-simple} for $\pqslope\in (3,4)$.



If on the other hand $\pqslope\in (2,3)$, the system has only rank $3$. Furthermore, using
the upper bound for the minimal seminorm does \emph{not} yet determine the system in this case.

But now recall from
Lemma~\ref{lem:symmetries-res}~(\ref{lem:sym-res-negate-p}) that
$\res_{p,q}\doteq \res_{-p+4q,q}$. Furthermore, $2<p/q<3$ if and only if
$1<-p/q+4<2$, and so the simplicity of the non-trivial roots of $\res_{p,q}$
when $\pqslope\in (2,3)$
follows from the simplicity of the non-trivial roots of $\res_{-p+4q,q}$ when
$\frac{-p+4q}{q}\in (1,2)$,
which we have already proven in section~\ref{sec:CS-positive-0-2-odd-filling}. We
can therefore already assume the validity of Corollary~\ref{cor:roots-of-res-are-simple} for
$\pqslope\in (2,3)$ and use that $s$ actually \emph{attains} the upper bound from
table~\ref{tab:number-of-conjug-classes-preps-Wpq}, i.e.\ $s=p+4q-3$.
Solving the linear system with this information, we then obtain the same
solution as given in equation \eqref{eq:CS-positive-3-4-odd-exact-solution},
completing the proof of Theorem~\ref{thm:total-cs-norm-Wpq} in this case.



\subsection{\texorpdfstring{Assuming $\pqslope\in (4,6)\cup (6,\infty)$,  $p$ odd}%
{Assuming 4<p/q<6 and 6<p/q, p odd}}
\label{sec:CS-positive-4-infty-odd-filling}

In this case, we obtain the following distances between the relevant slopes:
%
\begin{table}[h!]
  \centering
  \begin{tabular}{C|CCCCCCC}
    \Delta & \beta_{1}=4 & \beta_{2}=\frac{4q}{p-2q} & \beta_{3}=0 & 1 & 2 & 3 & \mu=\infty\\\hline
    \beta_{1} & 0 & 4p-12q & 4    & 3 & 2 & 1 & 1\\
    \beta_{2} &   &   0    & 4q   & |p-6q| & 2p-8q & 3p-10q & p-2q\\
    \beta_{3} &   &        & 0    & 1 & 2 & 3 & 1
  \end{tabular}
  \\[2ex]\caption{Distances between slopes in $\Wpq$ for $\pqslope\in
    (4,\infty)$, $p$ odd}
  \label{tab:distances-positive-4-infty-odd-filling}
\end{table}

Similar to the previous cases, we obtain the corresponding linear system as:
\begin{subequations}\label{eq:CS-positive-4-infty-odd}
  \begin{alignat}{2}
    s + 2(|p-6q|-1)  &= \CSnorm{1}      &&= 3a_{1} + |p-6q| a_{2}  + a_{3}\\
    s + 3(p-4q-1)    &= \CSnorm{2}      &&= 2a_{1} + (2p-8q)a_{2} + 2a_{3}\\
    s + 4(p-3q-1)    &= \CSnorm{3}      &&= a_{1}  + (3p-10q)a_{2}+ 3a_{3}\\
    s                &= \CSnorm{\infty} &&= a_{1}  + (p-2q)a_{2}  + a_{3}
  \end{alignat}
\end{subequations}
If $\pqslope\in (4,6)$, this system has rank $4$, and we can determine the
unique solution without further assumptions:\begin{equation}
  \label{eq:CS-positive-4-6-odd-exact-solution}
  \begin{pmatrix}
    a_{1}\\a_{2}\\a_{3}\\s
  \end{pmatrix}=
  \begin{pmatrix}
    p-2q-1\\
    2 \\
    2q-2\\
    3p-4q-3
  \end{pmatrix}.
\end{equation}
This then proves Theorem~\ref{thm:total-cs-norm-Wpq} and
Corollary~\ref{cor:roots-of-res-are-simple} for $\pqslope\in (4,6)$.

If on the other hand $\pqslope\in (6,\infty)$, the system has only rank
$3$. Furthermore, using the upper bound for the minimal seminorm does
\emph{not} yet determine the system in this case.

But the same argument as used in the previous subsection holds here: By
Lemma~\ref{lem:symmetries-res}~(\ref{lem:sym-res-negate-p}), we have that
$\res_{p,q}\doteq \res_{-p+4q,q}$. Furthermore, $p/q>6$ if and only if
$-p/q+4<-2$, and so the simplicity of the non-trivial roots of $\res_{p,q}$
when $\pqslope>6$
follows from the simplicity of the non-trivial roots of $\res_{-p+4q,q}$ when $\frac{-p+4q}{q}<-2$,
which we have already proven in section~\ref{sec:CS-negative-odd-filling}. We
can therefore already assume Corollary~\ref{cor:roots-of-res-are-simple} for
$\pqslope\in (6,\infty)$ and use that $s$ actually \emph{attains} the upper bound from
table~\ref{tab:number-of-conjug-classes-preps-Wpq}, i.e.\ $s=3p-4q-3$.
Solving the linear system with this information, we then obtain the same
solution as given in equation \eqref{eq:CS-positive-4-6-odd-exact-solution},
completing the proof of Theorem~\ref{thm:total-cs-norm-Wpq}.



\bibliographystyle{amsalpha}

\bibliography{references}

\end{document}